\documentclass{article}
\usepackage[dvips]{graphicx}
\usepackage{amsthm}  
\usepackage{amssymb}
\usepackage{amsmath}
\usepackage{enumerate, comment}
\usepackage{mathptmx}

\newcommand{\dt}{\circle*{3}}
\newcommand{\dnarr}{
\begin{picture}(5,15)
\put(3,15){\vector(0,-1){15}}
\end{picture}
}

\newcommand{\mrm}[1]{\mathrm{#1}}

\newcommand{\Aut}{\mrm{Aut}}

\newcommand{\integer}[1]{\lfloor #1 \rfloor}

\theoremstyle{definition}
\newtheorem{defn}{Definition}[section]
\newtheorem{rmk}[defn]{Remark}

\newtheorem{note}[defn]{Note}

\newtheorem{exa}[defn]{Example}
\theoremstyle{plain}
\newtheorem{thm}[defn]{Theorem}
\newtheorem*{thma}{Theorem A}
\newtheorem*{thmb}{Theorem B}

\newtheorem{lem}[defn]{Lemma}
\newtheorem{prop}[defn]{Proposition}

\begin{document}

\title{Cubic graphs with most automorphisms}
\author{Michael A. van Opstall, R\u{a}zvan Veliche}

\maketitle

\section{Introduction}

Let $G$ be a connected simple cubic graph; $\Aut~G$ denotes its 
automorphism group. Let $n$ be
half the number of vertices of $G$. We define the {\em arithmetic genus} of a 
(possibly disconnected) graph
 as $e-v+1$ where $e$ is the number of edges, $v$ the number of
vertices of $G$). For a connected simple 
cubic graph, $g=n+1$. The definition of arithmetic genus is motivated by
the following: to a projective nodal curve with rational components one may 
associate a so-called {\em dual graph}; the
arithmetic genus of the curve is the arithmetic genus of this graph. We 
abbreviate arithmetic genus to ``genus'' everywhere in this article, although
this is at variance with standard graph theory terminology. We trust that this
will not actually be confusing.

A bound on the order of $\Aut~G$ was obtained in \cite{wormald:79}, where it
is shown that $|\Aut~G|$ divides $3n\cdot 2^n$. However, it can be easily 
checked by consulting a list of cubic graphs\footnote{For example, Beresford's
Gallery of Cubic Graphs, which can be viewed at 
{\tt http://www.mathpuzzle.com/BeresfordCubic.html}.
This gallery is
complete for graphs of at most twelve vertices.}, that this bound is only 
rarely attained (in fact, it is only attained for graphs with four or six
vertices). Thus a natural problem is to find a sharp bound for the 
order of $\Aut~G$. We solve this problem by the following:

\begin{thm} Assume $g\geq 9$; set 
$l(g)=\min\{k| g=\sum\limits_{i=1}^k a_i\cdot 2^{n_i}, a_i\in \{1,3\}\}$, 
and set $o(g)=g-l(g)$.
\begin{itemize}
\item if $g=9\cdot 2^m+s$ ($s=0,1,2$) ($m\geq 0$) except $g=10,11,19,20,38$, 
then $|\Aut~G|\leq 3\cdot 2^{o(g)}$
\item if $g=3\cdot 2^m+s$ ($s=0,1,2$) ($m\geq 2$), or $g=9(2^m+2^p)$ (with 
$|m-p|\geq 5$) or if $g=10,11,19,20,38$, 
then $|\Aut~G|\leq \frac{3}{2}\cdot 2^{o(g)}$
\item if $g=5\cdot a\cdot 2^m+1$ (where $a=1$ or $3$, $m\geq 2$), then 
$|\Aut~G|\leq \frac{5}{4}\cdot 2^{o(g)}$ 
\item otherwise, $|\Aut~G|\leq 2^{o(g)}$
\end{itemize}
Moreover, these bounds are sharp; an explicit construction of graphs 
attaining the bounds in each case will be given in a subsequent section 
(see \ref{candidates}).
\end{thm}

The graphs with maximal automorphism group for $g\leq 8$ will be listed in
a table below.

This work was motivated by our earlier work \cite{msc} on maximal order 
automorphism
groups of stable curves. Aaron Bertram asked us if we could bound the
automorphism groups of stable curves with all rational components. This is
equivalent to the problem of finding the maximal order automorphism groups
of cubic multigraphs (with a slightly modified notion of graph automorphism).
Such a result may indeed be obtained along the lines of this article and is
pursued in \cite{msc2}.

The basic idea is as follows: once the genus is large enough (larger than 
eight), the graphs with the most automorphisms should be as nearly trees
as possible. Of course a tree cannot be trivalent, and all trees have
 genus zero. Subject to these restrictions, we need to attach 
``appendages'' of positive genus to trees in an optimal way. One sees that
this is easiest when the appendages have the smallest possible genus. 
Restricting to simple graphs forces the appendages to have genus at least
three, which in turn forces us to consider graphs slightly more general
than trees for the ``cores'' of our graphs. 
If we consider cubic multigraphs,
then there are appendages of genus two (a triangle with one edge doubled)
and genus one (a loop). The answers to Bertram's question are graphs with
loops as appendages. We will not pursue these questions on non-simple graphs
in this article.

\begin{figure}[ht]
\begin{center}
$\begin{array}{c@{\hspace{1.5in}}c}

\begin{picture}(30,30)
\put(5,15){\dt}
\put(5,15){\line(1,-1){10}}
\put(5,15){\line(1,1){10}}
\put(5,15){\line(1,0){20}}
\put(15,5){\dt}
\put(15,5){\line(1,1){10}}
\put(15,15){\dt}
\put(15,15){\line(0,1){10}}
\put(15,25){\dt}
\put(15,25){\line(1,-1){10}}
\put(25,15){\dt}
\end{picture}
 &
\begin{picture}(50,50)
\put(5,35){\dt}
\put(5,35){\line(1,1){10}}
\put(5,35){\line(1,-1){10}}
\put(5,35){\line(2,-3){20}}
\put(15,45){\dt}
\put(15,25){\dt}
\put(25,5){\dt}
\put(25,5){\line(2,3){20}}
\put(15,45){\line(1,-1){20}}
\put(15,45){\line(1,0){20}}
\put(15,25){\line(1,1){20}}
\put(15,25){\line(1,0){20}}
\put(35,45){\dt}
\put(35,25){\dt}
\put(35,45){\line(1,-1){10}}
\put(45,35){\dt}
\put(35,25){\line(1,1){10}}
\end{picture}
\\
\mbox{A pinched tetrahedron}
& 
\mbox{A pinched $K_{3,3}$}
\end{array}$
\end{center}
\caption{Pinching illustrated.}
\label{pinchfig}
\end{figure}
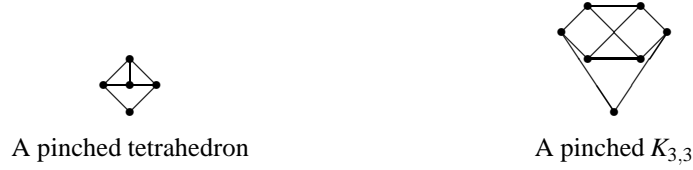

Our appendages are shown in Figure \ref{pinchfig}, and have genus three 
and four. The ``pinched tetrahedra'' were used previously, for example, in
an article of Wormald \cite{wormald:79}. A graph of genus 16 with 
$8^5=2^{o(16)}$ automorphisms may be constructed by attaching four copies
of the ``pinched $K_{3,3}$'' around the ends of a binary tree with
four leaves (with the root vertex removed so the graph is cubic). To reach
the bound given for genus 18, we factor 18 as $3\cdot 6$. Our goal is to
arrange six pinched tetrahedra around a core as symmetrically as possible. 
This is achieved by attaching them in three pairs to binary trees with
two leaves, and then arranging these binary trees around a star with four
vertices. 

The main idea is to construct a candidate graph in each genus whose
automorphism group has order equal to the bounds in the Main Theorem, thereby
giving a lower bound for the upper bound we seek. For a graph $G$ assumed to
be optimal, we
then remove the orbit of a well-chosen edge and attempt to proceed by
induction. However, removing this orbit in general causes problems - the
components of the remaining graph are not cubic, making them cubic may lead
to a non-simple graph, some components may be cycles, which cannot be made
cubic in any useful way. When one of these problems occurs, we will show that
it constrains the automorphism group of the graph so that its order is smaller
than that of the candidate graph's automorphism group. 

Graphs are assumed to be connected, unless otherwise stated. We will actually
have to work extensively with disconnected graphs, but this will be clearly
stated.

It is a pleasure to thank Marston Conder for his assistance in pointing out
the results of Goldschmidt's article, allowing us to significantly shorten
the elimination of edge-transitive graphs from consideration. Professor
Conder also made some helpful suggestions on an early version of the 
manuscript.

\section{Technicalities}

\begin{defn}
For a natural number $g$, define the functions
\begin{eqnarray*}
l(g)&=&\min\{k:g=\sum_{i=1}^k a_i\cdot 2^{n_i}, a_i\in\{1,3\}\} \\
o(g)&=&g-l(g).
\end{eqnarray*}
\end{defn}

The function $l$ may be computed as follows: expand $g$ in binary and starting
from the left, count the number of pairs $10$ or $11$, possibly adding one
to the total if after pairing there is a $1$ in the last digit. For example
$l(15)=2$, $l(21)=3$.

We will make extensive use of various inequalities involving the function 
$o$; we collect them here for convenience:

\begin{prop}\label{daily.inequalities}
\
\begin{itemize}
 \item $\frac{3}{2}k< 2^{o(k+1)}$ for all $k$.
 \item $k< 2^{2+o(\lceil \frac{k}{6}\rceil+1)}$ for all $k\geq 1$.
 \item $k<2^{o(\lceil \frac{2k}{3}\rceil +1)}$ for all $k$.
 \item $3k< 2^{o(2k-2)}$ for $k\geq 4$.
 \item $k\leq 2^{o(k-1)}$ for $k\geq 4$, with strict inequality for $k\geq 5$.
 \item $k\leq 2^{\lfloor \frac{k+1}{2}\rfloor}$, with strict inequality for 
$k\neq 2, 4$.
 \item $o(\lceil \frac{k}{2}\rceil +1)-\lfloor \frac{k}{6}\rfloor\geq o(\lceil 
\frac{k}{8}\rceil+1)$.
 \item $k\leq 2^{o(\lceil \frac{k}{2}\rceil+1)}$, with strict inequality
for $k\neq 2, 4, 8$.
\end{itemize}
\end{prop}

\begin{proof}
Straightforward.
\end{proof}

We list also some properties of the function $l$ which will be useful in what 
follows:

\begin{lem}\label{estimate.l.function}
We have the following (in)equalities:
\begin{itemize}
\item $l(a)=1$ if and only if $a=2^m$ or $a=3\cdot 2^m$ for some $m\geq 0$.
\item $l(a)\leq \lceil \frac{\log_2(s)}{2}\rceil$.
\item $l(a\cdot b)\leq 2l(a)\cdot l(b)$.
\item $l(2a)=l(a)$.
\item $o(2a)=o(a)+a$.
\item $l(3a)\leq 2l(a)$.
\item $l(a+1)\leq l(a)+1$.
\item $o(a)=1$ if and only if $a=2$; $o(a)\geq 2$ for $a\geq 3$.
\item $2^{l(a)}\leq 2\sqrt{a}$.
\item $l(a)=2$ and $l(3a)=4$ if and only if $a=3\cdot(2^m+2^p)$ with 
$m,p\geq 0$ and $|m-p|\geq 5$.
\end{itemize}
\end{lem}

\begin{proof}
Only the last two parts deserve some explanations: for $u$ such that 
$2^u\leq a <2^{u+1}$, the binary decomposition of $a$ will have $u+1$ digits; 
writing $a$ as sums in the definition of $l(a)$ effectively forms groups of 
at least two digits in this binary form, plus at most an extra one at the 
end; there are thus at most $\frac{u+2}{2}$ such groups, so 
$l(a)\leq \frac{u+1}{2}$; then $2^{l(a)}\leq 2\sqrt{2^u}\leq 2\sqrt{a}$. 
It is easy to see that the inequality is strict as soon as $a>1$.

For the last part, note that $l(a)=2$ means $a=b\cdot 2^m+c\cdot 2^m$ with 
$b,c\in\{1,3\}$. It is immediate to check that if at least one of $b,c$ is 
not $3$, then $l(3a)\geq 3$, and moreover, even if both $b=c=3$ one must have 
$|m-p|\geq 5$ for $l(3a)=4$ to happen.
\end{proof}

\begin{lem}\label{main.inequality}
The inequality $$sl(h)-l(sh)\geq \lfloor \frac{s+1}{2}\rfloor$$ is:
\begin{itemize}
 \item strict for any $h$ if $s=4,6$ or $s\geq 8$
 \item strict for $l(h)\geq 2$ and any $s\geq 2$, $s\neq 3$
 \item strict for $l(h)\geq 3$ and $s=3$
 \item equality for $l(h)=1$ and $s=2,5,7$, or $l(h)=2$ and $s=3$
 \item false for $s=1$ or $l(h)=1$ and $s=3$
\end{itemize}
\end{lem}
\begin{proof}
 We begin by noting that $l(sh)\leq 2l(s)l(h)$, so 
$sl(h)-l(sh)\geq l(h)(s-2l(s))=l(h)(2o(s)-s)$. Since $l(h)\geq 1$ we would 
like to see from what value of $s$ we have 
$2o(s)-s\geq \lfloor \frac{s+1}{2}\rfloor$, or equivalently 
$2o(s)\geq \lfloor \frac{3s+1}{2}\rfloor$.

If $2^u\leq s<2^{u+1}$ then 
$l(s)\leq \lceil \frac{u+1}{2}\rceil=\lfloor \frac{u}{2}\rfloor +1$. Then 
$2l(s)\leq 2\lfloor \frac{u}{2}\rfloor+2\leq u+2$, so 
$2o(s)=2s-2l(s)\geq 2s-u-2$. The inequality we would like to prove becomes 
$2s-u-2\geq \lfloor \frac{3s+1}{2}\rfloor $ or 
$s-\lfloor \frac{s+1}{2}\rfloor\geq u+2$. This in turn is implied by 
$\frac{s-1}{2}\geq u+2$ or $s\geq 2u+5$. Since $s\geq 2^u$ and $2^u>2u+5$ 
for $u\geq 4$ we see that the initial inequality is strict for 
$s\geq 16$.

Now one may construct a table of values for both sides of the inequality for 
values of $s$ up to $15$; we use the inequalities in 
(\ref{estimate.l.function}) above:

\begin{equation}\label{long.table}
\begin{array}{cccl}
s & sl(h)-l(sh) & \lfloor \frac{s+1}{2}\rfloor & \mrm{comments}\\
\hline
1 & 0 & 1 & \mrm{always\ false}\\
2 & l(h) & 1 & \mrm{not\ strict\ for\ }l(h)=1\\
3 & \geq l(h) & 2 & \mrm{strict\ for\ }l(h)\geq 3\\
4 & 3l(h) & 2 & \mrm{always\ strict}\\
5 & \geq 3l(h) & 3 & \mrm{not\ strict\ for\ }l(h)=1\\
6 & \geq 4l(h) & 3 & \mrm{always\ strict}\\
7 & \geq 4l(h) & 4 & \mrm{not\ strict\ for\ }l(h)=1\\
8 & 7l(h) & 4 & \mrm{always\ strict}\\
9 & \geq 7l(h) & 5 & \mrm{always\ strict}\\
10 & \geq 8l(h) & 5 & \mrm{always\ strict}\\
11 & \geq 8l(h) & 6 & \mrm{always\ strict}\\
12 & \geq 10l(h) & 6 & \mrm{always\ strict}\\
13 & \geq 10l(h) & 7 & \mrm{always\ strict}\\
14 & \geq 11l(h) & 7 & \mrm{always\ strict}\\
15 & \geq 11l(h) & 8 & \mrm{always\ strict}
\end{array}
\end{equation}

The lemma follows now immediately.
\end{proof}

In the future, we will denote by $A(s,h)$ the quantity 
$sl(h)-l(sh)-\lfloor \frac{s-1}{2}\rfloor$, and by $B(s,h)$ the quantity 
$sl(h)-l(sh)$; the lemma may be interpreted 
as giving ranges of $s$ and $h$ for which $A(s,h)\geq 1$.

\section{Eliminating edge-transitive graphs}

For an edge $e$ of a graph $G$, we will denote by $O(e)$ its orbit via the 
automorphism group of $G$. We use the word ``edge'' here to denote the graph
with two vertices joined by an edge, so that $O(e)$ is a graph. Define the
function $M(G)$ as the number of edges in a minimal orbit of an edge.

We refer to the star on four vertices simply as a ``star'', since we consider
no stars with more vertices.

\begin{lem}\label{structure.minimal.orbit}
Let $e$ be an edge of $G$ such that $O(e)$ has minimal order among all orbits 
of edges of $G$. Then only the following possibilities occur:
\begin{itemize}
\item $G=O(e)$;
\item $O(e)$ is a disjoint union of stars;
\item $O(e)$ is a disjoint union of edges;
\item $O(e)$ is a disjoint union of cycles; two such cycles are at distance 
at least two from each other.
\end{itemize}
\end{lem}

\begin{proof}
If is easy to see that if two stars in $O(e)$ have a common edge, then 
$G=O(e)$. Similarly, if two stars in $O(e)$ have a vertex in common, then 
either $G=O(e)$ or the third edge at that vertex will have an orbit of order 
smaller than that of $O(e)$ (which would be a contradiction to the choice of 
$e$).

Thus, if there is a star in $O(e)$, one of the first two possibilities occurs 
for $G$.

If no three edges in $O(e)$ share a common vertex, then either all edges in 
$O(e)$ are disjoint, or there are two edges $e_1$ (which may be assumed to be 
$e$, as $O(e)$ is acted upon transitively by $\Aut~G$) and $e_2$ in $O(e)$ 
with a common vertex $v$. Denote by $f$ the third edge of $G$ at $v$; $f$ is 
then not in $O(e)$. Denote by $w$ the other end of $e$.

If $v$ and $w$ are not in the same orbit of $\Aut~G$, then we see that 
$|O(e)|=2|O(v)|>|O(v)|\geq |O(f)|$, so we reach a contradiction to the 
choice of $e$. If however, $w\in O(v)$, then the existence of a cycle made 
of edges in $O(e)$ is immediate. Moreover, since $f\notin O(e)$, these cycles 
are disjoint.

Note that $|O(f)|\leq |O(e)|=|O(v)|$, with equality if and only if the ends 
of $f$ are not in the same orbit; in particular two cycles in $O(e)$ cannot 
be at distance one from each other (the edge between them, necessarily in the 
orbit of $f$, would have both endpoints in the same orbit).
\end{proof}

\begin{note}\label{e-f}
If the fourth situation above occurs, we will actually choose the edge $f$ 
and work with it in the arguments that follow; this is possible since $O(f)$ 
is also minimal, and may only be either a disjoint union of stars, or a 
disjoint union of edges. $f$ (or more precisely, its orbit) in this case will 
be called  {\em well-chosen}.
\end{note}

\begin{thm}\label{edge.transitive.bound}
An edge-transitive graph $G$ has at most $384(g-1)$ automorphisms.
\end{thm}

\begin{proof}
Tutte's papers \cite{tutte:47}\cite{tutte:59} on symmetric graphs give a 
bound of $48(g-1)$ for such graphs,
and Goldschmidt's results on semisymmetric graphs \cite{gold}
imply the bound in the statement of the theorem.
\end{proof}

We will use a couple of other functions frequently. Define 
\[
\mu(g)=\max \frac{|\Aut~G|}{2^{o(g)}},
\]
the maximum taken over all simple cubic graphs of genus $g$. For an edge
$e$ of a graph, define $\Aut_e' G$ to be the group of automorphisms
{\em preserving} (not necessarily fixing!) the edge $e$ -- that is, 
preserving the unordered pair of endpoints of $e$. Similarly, define
\[
\mu_1(g)=\max\frac{|\Aut_e' G|}{2^{o(g)}},
\]
the maximum here taken over all simple cubic graphs of genus $g$ and all edges
of these graphs. We are interested in the values of these functions for small
$g$ (when the optimal graphs are not our candidates). For a fixed graph $G$,
we define $\mu_1(G)$ similarly by taking the maximum over all edges of $G$.
Finally, set
\[
\pi(G)=\max|\Aut_e' G|,
\]
the maximum taken over all edges in $G$.

The following table may be compiled by inspection of lists of 
cubic graphs on a small number of vertices. 

\[\label{mu.table}
\begin{array}{lllll}
g & l(g) & o(g) & \mu(g) & \mu_1(g)\\
\hline
3 & 1 & 2 & 6 & 1 \\
4 & 1 & 3 & 9 & 1\\
5 & 2 & 3 & 6 & 1\\
6 & 1 & 5 & \frac{15}{4} & 1 \\
7 & 2 & 5 & 2 & 1\\
8 & 1 & 7 & \frac{21}{8} & 1
\end{array}
\]

\section{The candidate graphs}

We describe now the candidates $C_g$ for the cubic simple graphs with the
most automorphisms when $g\geq 9$. The definitions make sense for smaller
genus, but do not give the optimal graphs.

We need some non-standard terminology. If $G$ is an edge-transitive graph,
choose an edge $e$. Replace $e$ by two edges with one endpoint in common
and the other endpoints the former endpoints of $e$. We call this 
{\em pinching} $G$. This notion of pinching motivates the study of the function
$\mu_1$. If $G$ is not edge transitive, we must specify an edge when
pinching. 

A tree has a unique edge or vertex through which all
geodesics of maximal length pass; call this the {\em root}. If we attach a 
tree to
another graph at its root, and this root is an edge, we pinch the edge
and the new vertex is the point of attachment. If this attaching leads to an
a vertex of higher valence, we tacitly introduce an edge to correct the 
problem. In most cases, the meaning of ``attach'' is not confusing, since we
make the simplest attachment possible to keep the graph cubic. If there is
possibility of confusion, we will be more explicit.

\begin{defn}\label{candidates}
\

\begin{itemize}
\item Define the graphs $A_m$ as follows: attach a pinched tetrahedron to
each leaf of a binary tree with $2^m$ leaves. $A_m$ has genus 
$g=3\cdot 2^m$ ($m\geq 1$) and $2^{g-1}=2^{o(g)}$ automorphisms. For
$g=3\cdot 2^m$ ($m\geq 1$), define $C_g=A_m^\mrm{stab}$ (see the next section
for the definition of stabilization: in this case it ensures that the
graph is cubic by shrinking the central edge pair to an edge).
We note that $\pi(C_g)=2^{o(g)}$ in this case.

\item Define the graphs $B_m$ for $m\geq 2$ as follows: attach a pinched 
$K_{3,3}$ to each
leaf of a binary tree with $2^{m-2}$ leaves. $B_m$ has genus $g=2^m$ 
 and $2^{g-1}=2^{o(g)}$ automorphisms. For $g=2^m$ ($m\geq 3$), define
$C_g=B_m^\mrm{stab}$.
Note also that $\pi(C_g)=2^{o(g)}$ in this case.

\item If $g=3\cdot 2^m$ ($m\geq 2$), $C_g$ is defined by linking three copies 
of $B_m$ at 
their roots, by three edges, to a common root vertex. This is easily seen to 
yield a 
simple cubic graph with $\frac{3}{4}\cdot 2^g=\frac{3}{2}\cdot 2^{o(g)}$ 
automorphisms. If $g=3\cdot 2^m+1$, define $C_g$ by expanding the root of the 
previous tree 
into a triangle (with the $B_m$ attached at its vertices); this simple 
cubic graph has $\frac{3}{8}\cdot 2^g=\frac{3}{2}\cdot 2^{o(g)}$ 
automorphisms. Note that for this configuration, since $M(G) \geq 3$, we get 
$\pi(C_g)=\frac{1}{2} 2^{o(g)}$. 

\item If $g=9\cdot 2^m$, $C_g$ is defined by linking three copies of $A_m$ at 
their roots, by three edges, to a common root. This yields a simple cubic 
graph with 
$\frac{3}{4}\cdot 2^g=3\cdot 2^{o(g)}$ automorphisms. For $g=9\cdot 2^m+1$ 
we proceed as above by inserting a triangle at the root of the previous tree; 
again $|\Aut~C_g|=\frac{3}{8}\cdot 2^g=3\cdot 2^{o(g)} (m\geq 2)$ or 
$|\Aut~C_g|=\frac{3}{2}\cdot 2^{o(g)} (m=0,1)$. We note that 
$\pi(C_g)=2^{o(g)}$ for $m\geq 2$ and $\pi(C_g)=\frac{1}{2} 2^{o(g)}$ for 
$m=0,1$.

\item If $g=3\cdot 2^m+2$ ($m\geq 2$) or $g=9\cdot 2^m+2$ ($m\geq 0$), we 
obtain $C_g$ by 
attaching copies of $B_m$, respectively, $A_m$ to the valence two vertices of 
a $K_{2,3}$. In both 
cases the genus two $K_{2,3}$ at the core yields extra symmetry for a total of 
$\frac{3}{8}\cdot 2^g$ automorphisms. If $g=3\cdot 2^m+2$, this means 
$\frac{3}{2}\cdot 2^{o(g)}$ automorphisms, while if $g=9\cdot 2^m+2$ one gets 
$3\cdot 2^{o(g)}$ automorphisms for $m\geq 3$ and $\frac{3}{2}\cdot 2^{o(g)}$ 
automorphisms for $m=0,1,2$. As above, we note that $\pi(C_g)=2^{o(g)}$ for 
$m\geq 2$ and $\pi(C_g)=\frac{1}{2} 2^{o(g)}$ for $m=0,1,2$.

\item If $g=9(2^m+2^p)+s$ ($s=0,1,2$) with $|m-p|\geq 5$, $C_g$ is defined by
attaching $A_m$ to $A_k$ at their roots, and arranging three copies of this
configuration around a root which is a star, a triangle, or a $K_{2,3}$
depending on the value of $s$. This graph has $\frac322^{o(g)}$ automorphisms.
$\pi(C_g)=\frac122^{o(g)}$ in this case.

\item If $g=2^m+1$ ($m\geq 4$, as $m=3$ is covered above), we attach four 
copies of $B_m$ to the vertices of a square to obtain $C_g$; each of the 
quasi-trees 
will have genus $2^{m-2}$; the order of the automorphism group of this graph 
is then $8\cdot (2^{2^{m-2}-1})^4= 2^{g-2}=2^{o(g)}$. In this case 
$\pi(C_g)=\frac{1}{4}2^{o(g)}$.

\item If $g=5\cdot 2^m+1$ ($m\geq 2$) or $g=5\cdot 3\cdot 2^m+1$ ($m\geq 0$), 
then $C_g$ is a cycle of length five with copies of $B_m$, respectively $A_m$
at its vertices. In this case one gets $\frac{5}{4}\cdot 2^{o(g)}$ 
automorphisms.

\item $C_7$ is a pinched tetrahedron joined to a pinched $K_{3,3}$ by an
edge. 

\item In all other cases, in the binary decomposition of $g$ we have, counting
from the left, at least 
two pairs $11$ and/or $10$, plus a possible 1 left over. We look at the 
binary decomposition of $g$ and, from left to right, look for groups $11$ and 
$10$. We get a decomposition of $g$ into a sum 
of powers of two with coefficients one or three; for each part of $g$ of the 
form $2^m$ (with $m>1$) we take a $B_m$ of the corresponding genus, and 
for each part of $g$ of the form $3\cdot 2^m$ we take an $A_m$. 
If a 1 is 
left after this pairing, we replace the root of the last binary tree with a
triangle attached the the two branches, with a free edge attached to its
other branch. If a 10 is left, attach to the root of the last binary tree
an edge connected to a pair of triangles with a common side.

 In the end,  link each of these 
graphs to a distinct vertex of a path of length $l(g)-1$ using an edge. 
It is easy to 
see that the graph such obtained is cubic, simple, and of genus $g$; moreover, 
it is easy to compute that the order of the automorphism group of this graph 
is precisely $2^{o(g)}$.
\end{itemize}
\end{defn}

The final point in the definition will be called the {\em general case}, and
the others {\em exceptional}.

\begin{exa}
It is worth illustrating the general case with examples. Let $g=57$, so
the binary expansion of $g$ is $111001$. From left to right, there are two
pairs -- $11$ and $10$, and then a 1 is left over. We write 
$57=3\cdot 2^4+1\cdot 2^3+1\cdot 2^0$; $l(57)=3$, so $o(57)=54$. We are to take an $A_4$ 
and attach it to a
$B_3$, inserting a triangle in the middle. A simpler way to describe the
graph in this case is as follows: to two of the vertices of a triangle,
attach an $A_3$, and attach a $B_3$ to the third. 

$C_{56}$ is $C_{57}$ with the triangle collapsed to a point. There is no
change in automorphism group. $o(56)=o(57)=54$.

For $g=58$, we have
instead a pair of triangles sharing an edge, with an $A_4$ attached to the
free vertex of one, and the $B_3$ attached at the other free vertex.
$o(58)=55$, and there is an extra involution of $C_{58}$ coming from the
configuration of two triangles.
\end{exa}

\begin{rmk}\label{comparison}
\
\begin{itemize}
\item We note that the candidates above yield no more automorphisms than 
$3\cdot 2^{o(g)}$, and in most genera the bound is $2^{o(g)}$.  
\item We compare below the orders of the automorphism groups of the
candidates constructed above with the theoretical bound of $48(g-1)$ obtained 
by Tutte for symmetric graphs (there are no semisymmetric graphs of genus 
smaller than 28). For $9\leq g \leq 12$, methods similar to those of
Tutte give a bound of $24(g-1)$ for such graphs which will be used in
the table.

\[\label{low.genus.table}
\begin{array}{lllll}
g & o(g) & |\Aut~C_g| & \mrm{edge-transitive~bound} & \mrm{optimal} \\
\hline
3 & 2 & N/A & 96 & \mrm{tetrahedron~} (24)\\
4 & 3 & N/A & 144 & K_{3,3}~(72)\\
5 & 3 &  N/A & 192 & \mrm{cube~} (48)\\
6 & 5 & 32 & 240 & \mrm{Petersen~} (120)\\
7 & 5 & 32 & 288 & \mrm{no-name~} (64) \\
8 & 7 & 128 & 336 & \mrm{Heawood~} (336)\\
9 & 7 & {\bf 384} & 192 & \\
10 & 8 & {\bf 384} & 216 & \\
11 & 9 & {\bf 768} & 240 & \\
12 & 11 & {\bf 3072} &  264 & \\
13 & 11 & {\bf 3072} & 576 & \\
14 & 12 & {\bf 6144} & 624 & \\
15 & 13 & {\bf 8192} & 672 & \\
16 & 15 & {\bf 32768} & 720 &
\end{array}
\]
\end{itemize}
\end{rmk}

It is easy to see that for $g\geq 16$ we always have $2^{o(g)}> 384(g-1)$.

The conclusion so far: the graphs $C_g$ have more automorphisms than any 
edge-transitive graph as soon as $g\geq 9$.

Thus, we will be concerned in what follows with the two cases in which the 
minimal orbit of an edge is a disjoint union of stars or edges inside $G$.

As a consequence of (\ref{estimate.l.function}), we have the following:

\begin{prop}\label{growth.candidates}
$|\Aut~C_{g+1}|\geq |\Aut~C_g|$ except when $g=9\cdot 2^m+2$ ($m\geq 1$)
or $g=9(2^m+2^p)+2$ ($|m-p|\geq 5$). 
$|\Aut~C_{g+2}|>|\Aut~C_{g}|$ for any $g\geq 9$.
\end{prop}

\begin{proof}
Lemma \ref{estimate.l.function} shows that $o(g+1)\geq o(g)$.
\begin{itemize}
\item If $C_g$ has $2^{o(g)}$ automorphisms, then $C_{g+1}$ has at least 
$2^{o(g+1)}\geq 2^{o(g)}$.
\item If $C_g$ has $3\cdot 2^{o(g)}$ automorphisms, then $g=9\cdot 2^m+u$, 
where $u=0$ or $u=2$; in the second case $m\geq 3$ also. But then $C_{g+1}$, 
for $g+1=9\cdot 2^m+1$, has the same number of automorphisms as $C_g$, while 
for $g+1=9\cdot 2^m+3$, $C_{g+1}$ has $2^{o(g+1)}=2^{o(g)+1}$ automorphisms. 
Thus for $g=9\cdot 2^m+2$, $|\Aut~C_{g+1}|<|\Aut~C_g|$. Similar behavior 
occurs when $g=9(2^m+2^p)+2$ ($m$ and $p$ as in the hypotheses) with $p$ not
too small.
\item If $C_g$ has $\frac{3}{2} \cdot 2^{o(g)}$ automorphisms, then 
$o(g+1)\geq o(g)+1$ would show that $|\Aut~C_{g+1}| > |\Aut~C_g|$. We have:
\begin{itemize}
 \item If $g=3\cdot 2^m$, then $o(g+1)=o(g)$, but also 
$|\Aut~C_{g+1}|=\frac{3}{2}\cdot 2^{o(g+1)}$.
 \item If $g=3\cdot 2^m+u$ ($u=1,2$), then $o(g+1)=o(g)+1$.
 \item If $g=9\cdot 2^m+1$ ($n\geq 2$), then $o(g+1)=o(g)+1$.
 \item If $g=10,11,19$, in which case $g+1=11,12,20$ then again 
$o(g+1)\geq o(g)+1$ .
 \item If $g=20,38$, then as before we see that $|\Aut~C_{g+1}|<|\Aut~C_g|$.
 \end{itemize}
\end{itemize}

The previous calculations then show that $|\Aut~C_{g+1}|\geq|\Aut~C_g|$ 
except for the noted exceptions.

Since $o(g+2)\geq o(g)+1$, the second part follows immediately except in the 
case $g=9\cdot 2^m+2$; but a direct computation shows that $o(g+2)=o(g)+2$ 
for such $g$, so $|\Aut~C_{g+2}|\geq 4\cdot 2^{o(g)}>|\Aut~C_g|$ so we are 
done.
\end{proof}

\section{Some reductions}

\begin{defn}
We will call a simple cubic graph {\em optimal} if its automorphism group has 
maximal order among all simple cubic graphs of the same genus. We 
call such a graph {\em strictly optimal} if it is optimal and the minimal 
orbit of edges in it has minimal order among all optimal graphs of that 
genus.
\end{defn}

In this section we investigate the structure of a strictly optimal simple 
cubic graph of genus $g\geq 9$. The results of the previous 
sections show that $G$ cannot be edge-transitive. Consequently, we will pick 
a minimal orbit of an edge and try to understand its structure, and the 
structure it brings to the whole graph $G$.
 
A first step in this direction has been made by 
(\ref{structure.minimal.orbit}).

We will denote by $G'=G\setminus O(e)$ and by $g'$ the genus of 
$G'$. Due to the structure of $O(e)$, $G'$ has valence either two or three 
at each of its vertices. Consequently, $g'\geq 1$ and $g'=1$ if and only if 
$G'$ is a disjoint union of cycles.

\begin{rmk}\label{disc}

We will use a number of times a ``local surgery'' process, replacing subgraphs 
of the original graph $G$ with other subgraphs; the surgery is to be done 
throughout the orbit of the replaced subgraph. These graphs will have valence 
three at 
each of their vertices except at those that are 
in the same well-chosen orbit. We will reattach the replacements at the 
same points, to preserve regularity and avoid multiple edges. Each time we 
will keep track of the genus lost, and of the order of decrease/increase in 
the number of automorphisms. Since we are targeting only graphs that are 
optimal 
(maximal order of the automorphism subgroup), if by chance the orbit of the 
newly introduced subgraph is larger than that of the original subgraph, then 
we must have effectively/strictly increased the automorphism subgroup of the 
graph, so the estimates we use to show that the original one was not optimal 
still hold up. Thus we will assume tacitly that we are in the worst case 
scenario, where the new subgraph has orbit ``the same'' as the original, and 
argue usually by the number of elements in a minimal orbit of edges to reach 
a contradiction. More precisely, one can simply mark the vertices where the 
original graph was disconnected; the subgroup of automorphisms of the graph 
(after surgery) required to preserve the marking (at most permuting marked 
vertices among themselves) is the one we are really estimating. Most times 
this extra marking is not needed; in the few cases where it is, we will make 
it clear.
\end{rmk}

\begin{defn}
For a graph $G$ having valence at least two at each of its vertices and 
genus at least three we will denote by $G^\mrm{stab}$ the 
{\em stabilization} of $G$ (the terminology is motivated by algebraic
geometry). This is the graph obtained by replacing each maximal path, 
with interior vertices all of valence two by a simple edge with the same 
endpoints as the path. It is clear that $G^\mrm{stab}$ will have valence at 
least three at each of its vertices; however, it might have either loops or 
multiple edges.
\end{defn}

The possibility of loops and/or multiple edges inside $G^\mrm{stab}$ prohibits 
a direct induction; sharp bounds on such graphs may be obtained using the
methods of this article; they are greater than the bound for simple graphs
in every genus.

\begin{thma}
\label{main.theorem}
\
\begin{itemize}
\item If $G$ is a (strict optimal) simple cubic graph of genus $g\geq 9$, 
then $\mu(g)=\frac{|\Aut~G|}{2^{o(g)}}\leq 3$. Moreover, $\mu(g)\leq 1$ 
except when $g=a\cdot 2^m+b$, with 
$(a,b)\in \{(3,0),(3,1),(3,2),(9,0),(9,1),(9,2),(5,1),(15,1)\}$.

\item If $G$ is a (strict optimal) simple cubic graph of genus $g\geq 9$ 
with $|\Aut~G|>2^{o(g)}$, then $M(G)\geq 3$.

\item If $G$ is a simple cubic graph of genus $g\geq 9$, then 
$\pi(G)\leq 2^{o(g)}$ (or $\mu_1(g)\leq 1$).  
\end{itemize}
\end{thma}

Note that the third assertion of the theorem follows easily from the first 
two. Also note that Theorem A 
immediately implies the Main Theorem of the article. We will call Theorem A 
restricted to graphs of genus $g$ Theorem A$_g$.

The idea of the proof is to analyze the graph $G'$ left after removing a 
well-chosen minimal orbit $O(e)$ from $G$. $G'$ may be 
disconnected, $O(e)$ can be a disjoint union of stars or isolated edges, 
and $G'$ (or its components) might have loops or multiple edges after 
stabilization. Overall, $2^{o(g)}$ acts as a filter in each genus: 
if a graph does not have at least as many automorphisms, then it is not 
optimal, as (\ref{candidates}) shows. We will use estimates to determine 
precisely when (for what $g$, and for what structure of $O(e)$) a graph has
more than $2^{o(g)}$ automorphisms. Moreover, we will determine and use 
in the inductive process the order of the minimal orbit of edges (or at least 
some useful estimate). In what follows, $O(e)$ will {\em always} refer to a 
minimal orbit of disjoint stars or edges.

We call a vertex of $G'$ which has valence three a {\em stable} vertex. These 
might not exist. More precisely, these do not exist precisely in the 
components of $G'$ which are cycles. However, in this instance $G$ is easily 
seen not to be optimal; this will be shown inductively in Lemma 
\ref{all.cycles}.

In regard to the process of stabilization of the components of $G'$, the 
following lemmas will be useful:

\begin{lem}\label{unstable.path.length}
There can be at most two endpoints of edges in $O(e)$ on any unstable path 
with stable endpoints in $G'$.
\end{lem}

\begin{proof}
First note that if at least three such vertices exist on 
such a path, they cannot possibly be in the same orbit (there is in this case
a ``middle'' edge distinguished from the others). Since the endpoints of stars 
are in the same orbit, this means that we need to 
discuss only the possibility of isolated edges in $O(e)$ having three or more 
endpoints on a path with stable endpoints in $G'$. For such paths of length 
five or more (four or more contact points), it is clear that one of the edges 
in this path (the most central) will have an orbit of order less than that of 
$O(e)$, which is a contradiction. For a path of order four, the only 
possibility is that the middle contact point (vertex) is not in the orbit of 
the other two contact points. But then it has an orbit of order half of the 
other contact points, which is impossible, since both are endpoints of edges 
in $O(e)$.
\end{proof}

Moreover:

\begin{lem}\label{problems.stabilization}
Assume Theorem A$_h$ for all $9\leq h<g$.
Let $G$ be a strictly optimal simple cubic graph of genus $g$, and $O(e)$ a
minimal orbit of edges in $G$. Then no two edges in $O(e)$ may have endpoints 
on a path with stable endpoints in $G'=G\setminus O(e)$.
\end{lem}

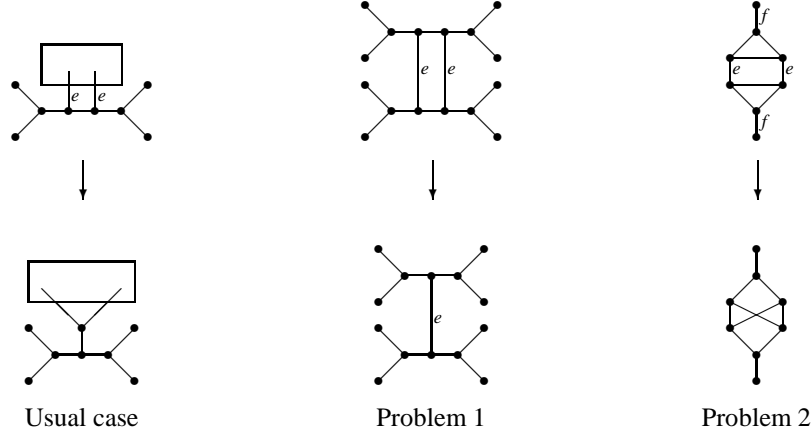
\begin{figure}[ht]
\begin{center}
$\begin{array}{c@{\hspace{1in}}c@{\hspace{1in}}c}
\begin{picture}(60,50)
\put(5,5){\dt}
\put(5,5){\line(1,1){10}}
\put(15,15){\dt}
\put(5,25){\dt}
\put(5,25){\line(1,-1){10}}
\put(15,15){\line(1,0){30}}
\put(25,15){\dt}
\put(35,15){\dt}
\put(45,15){\dt}
\put(45,15){\line(1,1){10}}
\put(45,15){\line(1,-1){10}}
\put(55,25){\dt}
\put(55,5){\dt}
\put(25,15){\line(0,1){15}}
\put(35,15){\line(0,1){15}}
\put(26,19){\makebox{\scriptsize $e$}}
\put(36,19){\makebox{\scriptsize $e$}}
\put(15,25){\framebox(30,15){}}
\end{picture}

&

\begin{picture}(60,70)
\put(5,5){\dt}
\put(5,5){\line(1,1){10}}
\put(15,15){\dt}
\put(5,25){\dt}
\put(5,25){\line(1,-1){10}}
\put(15,15){\line(1,0){30}}
\put(25,15){\dt}
\put(35,15){\dt}
\put(45,15){\dt}
\put(45,15){\line(1,1){10}}
\put(45,15){\line(1,-1){10}}
\put(55,25){\dt}
\put(55,5){\dt}

\put(25,15){\line(0,1){30}}
\put(35,15){\line(0,1){30}}
\put(26,29){\makebox{\scriptsize $e$}}
\put(36,29){\makebox{\scriptsize $e$}}

\put(5,35){\dt}
\put(5,35){\line(1,1){10}}
\put(15,45){\dt}
\put(5,55){\dt}
\put(5,55){\line(1,-1){10}}
\put(15,45){\line(1,0){30}}
\put(25,45){\dt}
\put(35,45){\dt}
\put(45,45){\dt}
\put(45,45){\line(1,1){10}}
\put(45,45){\line(1,-1){10}}
\put(55,55){\dt}
\put(55,35){\dt}
\end{picture}

& 
\begin{picture}(30,60)
\put(15,5){\dt}
\put(15,5){\line(0,1){10}}
\put(15,15){\dt}
\put(15,15){\line(-1,1){10}}
\put(15,15){\line(1,1){10}}
\put(5,25){\dt}
\put(5,35){\dt}
\put(5,25){\line(1,0){20}}
\put(5,25){\line(0,1){10}}
\put(25,25){\dt}
\put(25,25){\line(0,1){10}}
\put(25,35){\dt}
\put(5,35){\line(1,0){20}}
\put(5,35){\line(1,1){10}}
\put(15,45){\dt}
\put(15,45){\line(1,-1){10}}
\put(15,45){\line(0,1){10}}
\put(15,55){\dt}
\put(16,9){\makebox{\scriptsize $f$}}
\put(16,49){\makebox{\scriptsize $f$}}
\put(6,29){\makebox{\scriptsize $e$}}
\put(26,29){\makebox{\scriptsize $e$}}
\end{picture}

\\

\dnarr

&

\dnarr

&

\dnarr

\\

\begin{picture}(50,60)
\put(5,5){\dt}
\put(5,5){\line(1,1){10}}
\put(15,15){\dt}
\put(5,25){\dt}
\put(5,25){\line(1,-1){10}}
\put(15,15){\line(1,0){20}}
\put(25,15){\dt}
\put(35,15){\dt}
\put(25,15){\line(0,1){10}}
\put(25,25){\dt}
\put(25,25){\line(-1,1){15}}
\put(25,25){\line(1,1){15}}
\put(5,35){\framebox(40,15){}}
\put(35,15){\line(1,1){10}}
\put(35,15){\line(1,-1){10}}
\put(45,5){\dt}
\put(45,25){\dt}
\end{picture}

&

\begin{picture}(50,70)
\put(5,5){\dt}
\put(5,5){\line(1,1){10}}
\put(15,15){\dt}
\put(5,25){\dt}
\put(5,25){\line(1,-1){10}}
\put(15,15){\line(1,0){20}}
\put(25,15){\dt}
\put(35,15){\dt}
\put(35,15){\line(1,1){10}}
\put(35,15){\line(1,-1){10}}
\put(45,25){\dt}
\put(45,5){\dt}

\put(25,15){\line(0,1){30}}
\put(26,27){\makebox{\scriptsize $e$}}

\put(5,35){\dt}
\put(5,35){\line(1,1){10}}
\put(15,45){\dt}
\put(5,55){\dt}
\put(5,55){\line(1,-1){10}}
\put(15,45){\line(1,0){20}}
\put(25,45){\dt}
\put(35,45){\dt}
\put(35,45){\line(1,1){10}}
\put(35,45){\line(1,-1){10}}
\put(45,55){\dt}
\put(45,35){\dt}
\end{picture}

&

\begin{picture}(30,60)
\put(15,5){\dt}
\put(15,5){\line(0,1){10}}
\put(15,15){\dt}
\put(15,15){\line(-1,1){10}}
\put(15,15){\line(1,1){10}}
\put(5,25){\dt}
\put(5,35){\dt}
\put(5,25){\line(2,1){20}}
\put(5,25){\line(0,1){10}}
\put(25,25){\dt}
\put(25,25){\line(0,1){10}}
\put(25,35){\dt}
\put(5,35){\line(2,-1){20}}
\put(5,35){\line(1,1){10}}
\put(15,45){\dt}
\put(15,45){\line(1,-1){10}}
\put(15,45){\line(0,1){10}}
\put(15,55){\dt}
\end{picture}

\\

\mbox{Usual case} & \mbox{Problem 1} & \mbox{Problem 2}
\end{array}$
\end{center}
\caption{Stabilizations leading to a double edge. The problems are in the
first row, and their solutions in the second.}
\label{two.endpoints.plan}
\end{figure}

\begin{proof}
When two edges of $O(e)$ have contact points with a path with stable 
endpoints $p$ in $G'$, a naive (but effective) idea to try is to detach the 
two edges from the path, join their ``free'' ends to a common point, and join 
that common point to the middle of the original path (with the contact points 
removed) by a new edge $f$. This is illustrated in the first column of
Figure \ref{two.endpoints.plan}. 
The procedure should be carried on throughout $G$, 
in the orbit of the path $p$. In this way, the graph obtained has the same 
genus and it is easily seen to have at least as many automorphisms as the 
original one, and in fact the $O(f)$ will have order less than that of the 
$O(e)$ (this is easily seen to be true regardless of the structure of $O(e)$), 
so $G$ could not have been strictly optimal. There are only two possible 
cases when this procedure leads to double edges: 
\begin{enumerate}
\item {\bf Problem 1}: (does not happen when $O(e)$ is a disjoint union of 
stars): Two edges in $O(e)$ end up 
simultaneously with both ends on paths in the same orbit (see the middle
column of Figure \ref{two.endpoints.plan}); assume the path has distinct 
stable endpoints 
(for the other case, see {\bf Problem 2} below). However, in this instance 
replacing both edges by a single one yields a graph $\bar{G}$ which is simple, 
cubic, has no fewer automorphisms than $G$ (we do this on the whole orbit 
of the path at once) and has genus $g-\frac{k}{2}$ (where $k=|O(e)|$). 
Moreover, the orbit of the replacement edge has fewer members, contradicting 
the strict optimality of $G$. More precisely, (\ref{growth.candidates}) shows 
that, except in case $g=9\cdot 2^m+2$, $C_{g+1}$ has at least as many 
automorphisms as $C_g$. Then 
$|\Aut~G|\leq |\Aut(C_{g-\frac{k}{2}})|\leq |\Aut(C_g)|$; in case of equality 
throughout, the size of the minimal orbit helps establish the contradiction 
with the assumed strict optimality of $G$. This works except when $k=2$ and 
$g=9\cdot 2^m+3$, when the last inequality in the sequence above fails. 
However, in that case a direct argument can show that $G$ was not optimal; 
assuming $M(G)\geq 3$ for an optimal $G$ with $g=9\cdot 2^m+2$, and 
assuming that $|\Aut~H|\leq 2^{o(g)}$ for every non-optimal graph $H$ in the 
same genus, we see that $k=2$ and/or the orbit of $p$ having length two 
forces $G$ to be non-optimal, so 
$|\Aut~G|\leq 2^{o(g-1)}<2^{o(g)}=|\Aut~C_g|$.

\item {\bf Problem 2}: The procedure yields a double edge when the path $p$ 
has is a loop, starting and ending at the same stable point. Denoting by $f$ 
the third edge around this stable point, we see that the orbit of $f$ has 
order half that of $O(e)$ (contradiction!) unless the above-mentioned 
{\bf Problem 1} also occurs. But in this case the picture is as in the
third column of Figure \ref{two.endpoints.plan} and the whole configuration:  
may be replaced (see the same figure).
Since $f$ was not in the orbit of $e$, this stabilizes the graph locally, and 
in fact globally, as it is easy to see. Moreover, the automorphism group of 
the new graph (of the same genus as the original one) increases at 
least twofold, which is a contradiction to the optimality of $G$. Thus 
this problem does not occur in the optimal graphs.
\end{enumerate}
\end{proof}

Note that the previous lemma implies that any unstable path in $G'$ has
length two. Given this, we note, for future reference, what types of 
situations would lead to double (or triple) edges when stabilizing (the 
components of) 
$G'$. There are four classes, three labelled with roman numerals. A subscript 
on
a roman numeral indicates the length of a stable path between the stable 
endpoints
of an unstable path; if the subscript is a plus sign, the vertices are either 
connected 
by a stable path of length greater than two, or not connected by a stable 
path. The figures are drawn with the edges in $O(e)$ labelled ``$e$''.

\begin{figure}[ht]
\begin{center}
$\begin{array}{c@{\hspace{0.5in}}c@{\hspace{0.5in}}c@{\hspace{0.5in}}c}
\begin{picture}(50,30)
\put(5,15){\dt}
\put(5,15){\line(1,0){10}}
\put(15,15){\dt}
\put(15,15){\line(1,1){10}}
\put(15,15){\line(1,-1){10}}
\put(25,5){\dt}
\put(25,25){\dt}
\put(25,5){\line(0,1){20}}
\put(25,5){\line(1,1){10}}
\put(25,25){\line(1,-1){10}}
\put(35,15){\dt}
\put(35,15){\line(1,0){10}}
\put(45,15){\dt}
\put(9,16){\makebox{\scriptsize $e$}}
\put(39,16){\makebox{\scriptsize $e$}}
\end{picture}
&
\begin{picture}(80,30)
\put(5,15){\dt}
\put(5,15){\line(1,0){10}}
\put(15,15){\dt}
\put(15,15){\line(1,1){10}}
\put(15,15){\line(1,-1){10}}
\put(25,5){\dt}
\put(25,25){\dt}
\put(25,25){\line(1,-1){10}}
\put(25,25){\line(4,-1){40}}
\put(25,5){\line(1,1){10}}
\put(25,5){\line(4,1){40}}
\put(35,15){\dt}
\put(35,15){\line(1,0){10}}
\put(45,15){\dt}
\put(65,15){\dt}
\put(65,15){\line(1,0){10}}
\put(75,15){\dt}
\put(9,16){\makebox{\scriptsize $e$}}
\put(39,16){\makebox{\scriptsize $e$}}
\put(69,16){\makebox{\scriptsize $f$}}
\end{picture}
&
\begin{picture}(50,30)
\put(5,15){\dt}
\put(5,15){\line(1,0){10}}
\put(15,15){\dt}
\put(15,15){\line(1,1){10}}
\put(15,15){\line(1,-1){10}}
\put(25,5){\dt}
\put(25,25){\dt}
\put(25,5){\line(0,-1){5}}
\put(25,25){\line(0,1){5}}
\put(25,5){\line(1,1){10}}
\put(25,25){\line(1,-1){10}}
\put(35,15){\dt}
\put(35,15){\line(1,0){10}}
\put(45,15){\dt}
\put(9,16){\makebox{\scriptsize $e$}}
\put(39,16){\makebox{\scriptsize $e$}}
\end{picture}
&
\begin{picture}(30,30)
\put(15,5){\dt}
\put(15,5){\line(-1,1){10}}
\put(15,5){\line(0,1){20}}
\put(15,5){\line(1,1){10}}
\put(15,15){\dt}
\put(15,25){\dt}
\put(15,25){\line(-1,-1){10}}
\put(15,25){\line(1,-1){10}}
\put(5,15){\dt}
\put(25,15){\dt}
\end{picture}
\\
\mbox{I$_1$} & \mbox{I$_2$} & \mbox{I$_+$} & \mbox{$K_{2,3}$} 
\\
\begin{picture}(30,50)
\put(15,5){\dt}
\put(15,5){\line(0,1){10}}
\put(15,15){\dt}
\put(15,15){\line(-1,2){10}}
\put(15,15){\line(1,2){10}}
\put(15,25){\dt}
\put(15,25){\line(-1,1){10}}
\put(15,25){\line(1,1){10}}
\put(5,35){\dt}
\put(25,35){\dt}
\put(15,25){\line(0,1){20}}
\put(15,45){\dt}
\put(5,35){\line(1,1){10}}
\put(25,35){\line(-1,1){10}}
\put(16,33){\makebox{\scriptsize $e$}}
\put(16,9){\makebox{\scriptsize $f$}}
\end{picture}
&
\begin{picture}(50,30)
\put(5,15){\dt}
\put(5,15){\line(1,0){10}}
\put(15,15){\dt}
\put(15,15){\line(1,1){10}}
\put(15,15){\line(1,-1){10}}
\put(25,5){\dt}
\put(25,25){\dt}
\put(25,5){\line(0,1){20}}
\put(25,5){\line(1,1){10}}
\put(25,25){\line(1,-1){10}}
\put(35,15){\dt}
\put(35,15){\line(1,0){10}}
\put(45,15){\dt}
\put(26,14){\makebox{\scriptsize $e$}}
\end{picture}
&
\begin{picture}(50,30)
\put(5,15){\dt}
\put(5,15){\line(1,0){10}}
\put(15,15){\dt}
\put(15,15){\line(1,1){10}}
\put(15,15){\line(1,-1){10}}
\put(25,5){\dt}
\put(25,25){\dt}
\put(25,5){\line(0,1){20}}
\put(25,5){\line(1,1){10}}
\put(25,25){\line(1,-1){10}}
\put(35,15){\dt}
\put(35,15){\line(1,0){10}}
\put(45,15){\dt}
\put(9,16){\makebox{\scriptsize $e$}}
\end{picture}
&
\begin{picture}(50,30)
\put(5,5){\dt}
\put(5,5){\line(1,0){10}}
\put(5,25){\dt}
\put(5,25){\line(1,0){10}}
\put(15,5){\dt}
\put(15,25){\dt}
\put(15,5){\line(0,1){20}}
\put(15,5){\line(1,1){10}}
\put(15,25){\line(1,-1){10}}
\put(25,15){\dt}
\put(25,15){\line(1,0){10}}
\put(29,16){\makebox{\scriptsize $e$}}
\put(35,15){\dt}
\end{picture}
\\
\mbox{II$_2$} & \mbox{II$_+$} & \mbox{III$_2$} & \mbox{III$_+$} 
\end{array}$
\end{center}
\caption{The subgraphs that lead to multiple edges and loops.}
\label{nasty.subgraphs}
\end{figure}
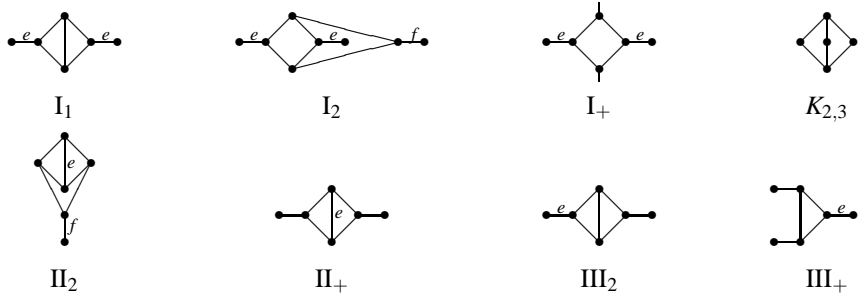

We will discuss when can these structures occur in a strictly optimal graph 
$G$. All the constructions are assumed to be done throughout the orbit of the 
unstable path simultaneously. We denote by $\bar{G}$ the graph obtained as a 
result of the surgery. The following definition is convenient.

\begin{defn}
Let $H$ be a graph which is cubic except for two vertices of valence two.
A {\em pseudocycle} is a cubic graph obtained by replacing the vertices in a 
cycle with copies of $H$.
\end{defn}

\begin{lem}
An optimal graph does not have a subgraph of type I$_1$.
\end{lem}

\begin{proof}
I$_1$ may only occur when the whole $G$ is a pseudocycle
 (otherwise the middle edge would have a 
shorter orbit). Direct computation shows that these graphs are not optimal for 
$g\geq 9$.
\end{proof}

\begin{lem}\label{situation.V}
A strictly optimal graph has no subgraphs of type I$_2$.
\end{lem}

\begin{proof}
{\bf Case 1:} A strictly optimal $G$ with $O(e)$ a disjoint union of stars
cannot have an I$_2$: it is clear that the orbit of the edge labelled $f$
in Figure \ref{nasty.subgraphs}  will have an orbit with fewer elements than 
$O(e)$.

{\bf Case 2:} When $O(e)$ is a disjoint union of $k\geq 2$ edges, the only
possible way that an I$_2$ could occur in $G$ is if members
of $O(e)$ alternate with I$_2$ to form a pseudocycle (otherwise again the 
edge labelled $f$ would have a smaller orbit). There are two subcases:
  \begin{enumerate}[a]

   \item The pseudocycles have exactly two I$_2$. We will create a new
graph whose minimal orbit is smaller, contradicting strict optimality
of $G$. Each pseudocycle is connected to the rest of $G$ by the edges labelled
$f$ in the figure. If these two edges actually coincide, the graph $G$ is
one pseudocycle, and direct computation excludes it from optimality. So we
may assume that each pseudocycle has two distinct edges connecting it to the
rest of the graph. A pseudocycle has local genus
five and local automorphism group of order eight (sixteen, if the two
$f$s may be flipped). Remove the pseudocycle and replace it with a triangle 
connected
by two vertices to the two $f$s and whose third vertex is connected to a
graph $A_1$. This does not change the genus or number of automorphisms of 
the graph, but the edge connecting the triangle to the $A_1$ now has a smaller 
orbit than the previous $e$.

   \item If the pseudocycle has at least three I$_2$, we replace it with a
cycle whose length is equal to the number of I$_2$. In this way the new 
graph, cubic and simple, has lost precisely $2^k$ 
automorphisms (permuting the stable endpoints of the unstable paths), but 
also lost $2k$ in genus. By induction, if $\bar{g}\geq 9$ we have 
$|\Aut~\bar{G}|\leq 3\cdot 2^{o(g-2k)}$ so 
$|\Aut~G|\leq 3\cdot 2^{k+o(g-2k)}\leq 3\cdot 2^{o(g)-o(k)}$. Since $k\geq 3$, 
the last quantity is clearly less than $2^{o(g)}$ so $G$ could not have been 
optimal. If however $\bar{g}\leq 8$, we must increase $\mu(G)$ to nine, so we 
can only derive the contradiction to the optimality of $G$ when $k\geq 6$; if 
$\mu(G)\leq 6$, then the contradiction happens as soon as $k\geq 4$. So we 
are left with a few cases to consider: 
  \begin{itemize}
   \item $\mu(G)=9$ and $k=4,5$. Then $\bar{g}=4$, so $g=12$ or $14$; 
$|\Aut~\bar{G}|\leq 72$. In the first case, 
$|\Aut~G|\leq 72\cdot 2^4<6\cdot 2^9=|\Aut~C_{12}|$ so $G$ was not optimal; 
in the second case $|\Aut~G|\leq 72\cdot 2^5<2^{12}=|\Aut~C_{14}|$ so again 
$G$ was not optimal
   \item $\mu(G)\leq 6$ and $k=3$. Then $g\leq 14$, and the table 
(\ref{low.genus.table}) shows that $|\Aut~G|<|\Aut~C_g|$ so $G$ is not optimal.
  \end{itemize}
  \end{enumerate}
\end{proof}

\begin{lem}
Assume Theorem $A_h$ for $9\leq h<g$. Then
an optimal graph of genus $g$ does not have a subgraph of type I$_+$.
\end{lem}

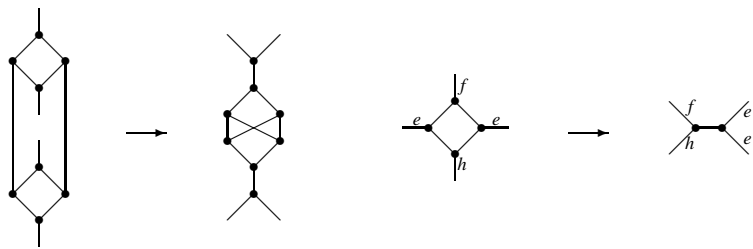
\begin{figure}[ht]
\begin{center}
$\begin{array}{c@{\hspace{0.25in}}c@{\hspace{0.25in}}c@{\hspace{0.5in}}c@{\hspace{0.25in}}c@{\hspace{0.25in}}c}
\begin{picture}(30,100)
\put(15,5){\line(0,1){10}}
\put(15,15){\dt}
\put(15,15){\line(-1,1){10}}
\put(15,15){\line(1,1){10}}
\put(5,25){\dt}
\put(25,25){\dt}
\put(5,25){\line(1,1){10}}
\put(25,25){\line(-1,1){10}}
\put(15,35){\dt}
\put(15,35){\line(0,1){10}}
\put(5,25){\line(0,1){50}}
\put(25,25){\line(0,1){50}}
\put(5,75){\dt}
\put(25,75){\dt}
\put(5,75){\line(1,1){10}}
\put(5,75){\line(1,-1){10}}
\put(25,75){\line(-1,1){10}}
\put(25,75){\line(-1,-1){10}}
\put(15,65){\dt}
\put(15,65){\line(0,-1){10}}
\put(15,85){\dt}
\put(15,85){\line(0,1){10}}
\end{picture}
&
\begin{picture}(15,5)(0,-45)
\put(0,3){\vector(1,0){15}}
\end{picture}
&
\begin{picture}(30,80)(0,-10)
\put(5,5){\line(1,1){10}}
\put(25,5){\line(-1,1){10}}
\put(15,15){\dt}
\put(15,15){\line(0,1){10}}
\put(15,25){\dt}
\put(15,25){\line(-1,1){10}}
\put(15,25){\line(1,1){10}}
\put(5,35){\dt}
\put(25,35){\dt}
\put(5,35){\line(0,1){10}}
\put(5,35){\line(2,1){20}}
\put(25,35){\line(0,1){10}}
\put(25,35){\line(-2,1){20}}
\put(5,45){\dt}
\put(25,45){\dt}
\put(5,45){\line(1,1){10}}
\put(25,45){\line(-1,1){10}}
\put(15,55){\dt}
\put(15,55){\line(0,1){10}}
\put(15,65){\dt}
\put(15,65){\line(-1,1){10}}
\put(15,65){\line(1,1){10}}
\end{picture}
&
\begin{picture}(50,50)(0,-25)
\put(5,25){\line(1,0){10}}
\put(15,25){\dt}
\put(15,25){\line(1,-1){10}}
\put(15,25){\line(1,1){10}}
\put(25,35){\dt}
\put(25,15){\dt}
\put(25,35){\line(0,1){10}}
\put(25,15){\line(0,-1){10}}
\put(25,35){\line(1,-1){10}}
\put(25,15){\line(1,1){10}}
\put(35,25){\dt}
\put(35,25){\line(1,0){10}}
\put(9,26){\makebox{\scriptsize $e$}}
\put(39,26){\makebox{\scriptsize $e$}}
\put(26,9){\makebox{\scriptsize $h$}}
\put(26,39){\makebox{\scriptsize $f$}}
\end{picture}
&
\begin{picture}(15,5)(0,-45)
\put(0,3){\vector(1,0){15}}
\end{picture}
&
\begin{picture}(40,30)(0,-35)
\put(5,5){\line(1,1){10}}
\put(15,15){\dt}
\put(5,25){\line(1,-1){10}}
\put(15,15){\line(1,0){10}}
\put(25,15){\dt}
\put(25,15){\line(1,1){10}}
\put(25,15){\line(1,-1){10}}
\put(11,7){\makebox{\scriptsize $h$}}
\put(11,21){\makebox{\scriptsize $f$}}
\put(33,9){\makebox{\scriptsize $e$}}
\put(33,19){\makebox{\scriptsize $e$}}
\end{picture}
\end{array}$
\end{center}
\caption{Surgeries for graphs of type I$_+$.}
\label{iplusnotsquare}
\end{figure}

\begin{proof}
{\bf Case 1:} Suppose $O(e)$ is a disjoint union of $k$ stars, so the 
components of $G'$ are all isomorphic. If the components of $G'$ pairs
of squares connected by two edges, diagonally opposite each other, we may
replace these components as in Figure \ref{iplusnotsquare} to gain 
automorphisms without changing the genus, contradicting optimality.

Otherwise, one may collapse each I$_+$ as in Figure \ref{iplusnotsquare}.
The resulting graph is cubic simple, and has at least as many automorphisms 
as $G$; however, its genus is $g-\frac{3k}{2}$, and then 
(\ref{growth.candidates}) and (\ref{mu.table}) show that $G$ could not have 
been optimal.

{\bf Case 2:} The same surgery as above may be done when $O(e)$ is a disjoint 
union of $k$ isolated edges (under the same restriction as above).
  \begin{itemize}
   \item If the two unstable vertices of an I$_+$ are not in the same orbit, 
then the existence of a pseudocycle formed of edges in $O(e)$ and I$_+$ 
to which they are incident is immediate; in fact the whole orbit $O(e)$ will 
be partitioned in edges arranged in such isomorphic pseudocycles. If a 
pseudocycle contains at least three I$_+$, do the same surgery as above; the 
decrease in genus overall is precisely $k\geq 2$, so again $G$ could not have 
been optimal. If a pseudocycle contains only two I$_+$, then the surgery is 
modified to unite by an edge the two vertices to which the eyes were 
contracted. This time the decrease in genus is at least three, so again $G$ 
could not have been optimal.

   \item If the two unstable vertices of an eye are in the same orbit, but 
the other ends of the two incident edges are not in their orbit, then the 
surgery above can be done at the orbit of the unstable path without leading 
to double edges; the decrease in genus is  $\frac{k}{2}$ (it is easy to see 
that $k$ must be even), and the only case requiring consideration is when 
$k=2$ (when the difference in genus is only one) and $\mu(\bar{G})>1$; but 
then $M(\bar{G})\geq 3$ by induction, so we reach a contradiction.

   \item If the ends of the edges in $O(e)$ are in the same orbit, we see 
again the existence of pseudocycles in $O(e)$ and we continue the argument as 
above; $G$ could not be optimal.
  \end{itemize}

In the case of ``cylinders'' (as in Case 1, the left side of the figure) we
may replace as in the case when $O(e)$ was a disjoint union of stars, and see
that $G$ was not optimal. This completes the proof.
\end{proof}

\begin{lem}
An optimal contains at most one subgraph of type III$_+$. 
Furthermore, if a graph of type III$_2$ occurs in an optimal graph, it is
unique, so may be considered as a subgraph II$_+$.
\end{lem}

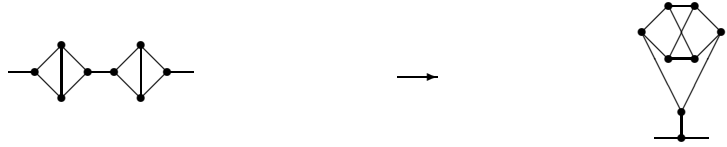
\begin{figure}[ht]
\begin{center}
$\begin{array}{c@{\hspace{1in}}c@{\hspace{1in}}c}
\begin{picture}(80,30)(0,-15)
\put(5,15){\line(1,0){10}}
\put(15,15){\dt}
\put(15,15){\line(1,1){10}}
\put(15,15){\line(1,-1){10}}
\put(25,5){\dt}
\put(25,25){\dt}
\put(25,5){\line(0,1){20}}
\put(25,5){\line(1,1){10}}
\put(25,25){\line(1,-1){10}}
\put(35,15){\dt}
\put(35,15){\line(1,0){10}}
\put(45,15){\dt}

\put(45,15){\line(1,1){10}}
\put(45,15){\line(1,-1){10}}
\put(55,5){\dt}
\put(55,25){\dt}
\put(55,5){\line(0,1){20}}
\put(55,5){\line(1,1){10}}
\put(55,25){\line(1,-1){10}}
\put(65,15){\dt}
\put(65,15){\line(1,0){10}}
\end{picture}

&
\begin{picture}(15,5)(0,-25)
\put(0,3){\vector(1,0){15}}
\end{picture}
&

\begin{picture}(40,60)
\put(10,5){\line(1,0){20}}
\put(20,5){\dt}
\put(20,5){\line(0,1){10}}
\put(20,15){\dt}
\put(20,15){\line(-1,2){15}}
\put(20,15){\line(1,2){15}}
\put(5,45){\dt}
\put(5,45){\line(1,1){10}}
\put(5,45){\line(1,-1){10}}
\put(15,55){\dt}
\put(15,35){\dt}
\put(15,55){\line(1,0){10}}
\put(15,35){\line(1,0){10}}
\put(25,55){\dt}
\put(25,35){\dt}
\put(15,35){\line(1,2){10}}
\put(15,55){\line(1,-2){10}}
\put(25,35){\line(1,1){10}}
\put(25,55){\line(1,-1){10}}
\put(35,45){\dt}
\end{picture}
\end{array}$
\end{center}
\caption{Surgery for graphs of type III$_2$.}
\label{breaking.iii2}
\end{figure}

\begin{proof}
The argument for subgraphs of type III$_+$ follows the line of those above. 
In this case, the surgery is to collapse the triangle in each III$_+$ to a
point. The resulting graph has the same automorphisms, and lower genus, which
contradicts optimality as above if there were multiple III$_+$ in the graph (if
there is only one, the genus does not decrease enough to apply 
(\ref{growth.candidates})). Note that these subgraphs do appear in certain
of the $C_g$.

The middle edge of the graph III$_2$ moves as much as the edge labelled $e$. 
If the orders of their orbits are equal, we may shift attention to the middle
edge and think of the III$_2$ as a II$_+$. In this case, the result follows
from the next lemma. 

The only way that the orbit of $e$ could be smaller than the orbit of the
middle edge is if a configuration as in Figure \ref{breaking.iii2} occurs, and 
the same figure gives a surgical solution to this problem.
\end{proof}

\begin{lem} 
An optimal graph may have at most one II$_+$.
\end{lem} 

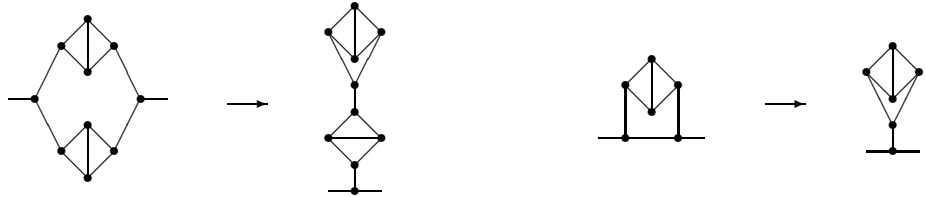
\begin{figure}[ht]
\begin{center}
$\begin{array}{c@{\hspace{0.25in}}c@{\hspace{0.25in}}c@{\hspace{1in}}c@{\hspace{0.25in}}c@{\hspace{0.25in}}c}
\begin{picture}(70,70)(0,-5)
\put(5,35){\line(1,0){10}}
\put(15,35){\dt}
\put(15,35){\line(1,2){10}}
\put(15,35){\line(1,-2){10}}
\put(25,55){\dt}
\put(25,15){\dt}
\put(25,55){\line(1,1){10}}
\put(25,55){\line(1,-1){10}}
\put(25,15){\line(1,1){10}}
\put(25,15){\line(1,-1){10}}
\put(35,65){\dt}
\put(35,45){\dt}
\put(35,5){\dt}
\put(35,25){\dt}
\put(35,65){\line(1,-1){10}}
\put(35,65){\line(0,-1){20}}
\put(35,45){\line(1,1){10}}
\put(35,25){\line(1,-1){10}}
\put(35,25){\line(0,-1){20}}
\put(35,5){\line(1,1){10}}
\put(45,55){\dt}
\put(45,15){\dt}
\put(45,55){\line(1,-2){10}}
\put(45,15){\line(1,2){10}}
\put(55,35){\dt}
\put(55,35){\line(1,0){10}}
\end{picture}
&
\begin{picture}(15,5)(0,-35)
\put(0,3){\vector(1,0){15}}
\end{picture}
&
\begin{picture}(30,80)
\put(5,5){\line(1,0){20}}
\put(15,5){\dt}
\put(15,5){\line(0,1){10}}
\put(15,15){\dt}
\put(15,15){\line(1,1){10}}
\put(15,15){\line(-1,1){10}}
\put(5,25){\dt}
\put(5,25){\line(1,0){20}}
\put(25,25){\dt}
\put(5,25){\line(1,1){10}}
\put(25,25){\line(-1,1){10}}
\put(15,35){\dt}
\put(15,35){\line(0,1){10}}
\put(15,45){\dt}
\put(15,45){\line(-1,2){10}}
\put(15,45){\line(1,2){10}}
\put(5,65){\dt}
\put(25,65){\dt}
\put(5,65){\line(1,1){10}}
\put(5,65){\line(1,-1){10}}
\put(15,75){\line(0,-1){20}}
\put(15,75){\dt}
\put(15,75){\line(1,-1){10}}
\put(15,55){\dt}
\put(15,55){\line(1,1){10}}
\end{picture}
&
\begin{picture}(50,40)(0,-20)
\put(5,5){\line(1,0){40}}
\put(15,5){\dt}
\put(15,5){\line(0,1){20}}
\put(35,5){\dt}
\put(35,5){\line(0,1){20}}
\put(15,25){\dt}
\put(35,25){\dt}
\put(15,25){\line(1,1){10}}
\put(15,25){\line(1,-1){10}}
\put(25,35){\dt}
\put(25,15){\dt}
\put(25,35){\line(1,-1){10}}
\put(25,15){\line(1,1){10}}
\put(25,15){\line(0,1){20}}
\end{picture}
&
\begin{picture}(15,5)(0,-35)
\put(0,3){\vector(1,0){15}}
\end{picture}
&
\begin{picture}(30,50)(0,-15)
\put(5,5){\line(1,0){20}}
\put(15,5){\dt}
\put(15,5){\line(0,1){10}}
\put(15,15){\dt}
\put(15,15){\line(-1,2){10}}
\put(15,15){\line(1,2){10}}
\put(5,35){\dt}
\put(25,35){\dt}
\put(5,35){\line(1,1){10}}
\put(5,35){\line(1,-1){10}}
\put(15,45){\line(0,-1){20}}
\put(15,45){\dt}
\put(15,45){\line(1,-1){10}}
\put(15,25){\dt}
\put(15,25){\line(1,1){10}}
\end{picture}
\end{array}$
\end{center}
\caption{Surgeries for multiple II$_2$.}
\label{no.two.theta}
\end{figure}

\begin{proof}
First of all, replacing each II$_+$ by a single edge cannot produce a double 
edge or a loop (producing a loop would mean that we had a II$_2$). This is 
because either two such II$_+$ share the vertices at distance one from their 
stable ends, or there is a shortcut (edge) between those vertices (at 
distance one from their stable ends); we note that, if three such II$_+$ 
share the vertices of distance one from their stable endpoints, this 
configuration is the whole graph, of genus eight, and excluded from 
consideration. The surgeries in these two cases are depicted in Figure
\ref{no.two.theta}. In the first case, the order of the orbit of a minimal
edge is decreased without changing genus or automorphisms, contradicting
strict optimality. In the second case, the surgery produces a graph with 
a larger automorphism group.

Then, if $k=|O(e)|\geq 2$, the II$_+$ could not have adjacent vertices: more 
than three II$_+$ could only be adjacent if they form a cycle (the whole of 
$G$!) due to the requirement that their middle edges should be in the same 
orbit; then $|\Aut~G|\leq 2k\cdot 2^k<2^{o(2k+1)}$ for $k\geq 5$ and 
$|\Aut~G|<|\Aut~C_9|$ for $k=4$. And if only two II$_+$ would be adjacent, 
the same surgery as in the previous lemma (see Figure \ref{breaking.iii2}) 
can be done.

Thus if more than two II$_+$ exist in a strictly optimal $G$, they are not 
adjacent, and the vertices at distance one from their stable endpoints are 
not neighbors; replace then the II$_+$ by single edges; the automorphism 
group decreased in order by a factor of $2^k$, the resulting graph is simple 
and cubic, and of genus $g-2k$. Then the same discussion as in the case of 
the I$_2$ shows that $G$ could not have been optimal.
\end{proof}

To eliminate most subgraphs of type II$_2$ from consideration, we refine
our notion of well-chosen. Note that our candidate graphs contain many copies 
of II$_2$, but with $e$ not minimal. We will avoid the multiple edges
which occur as a result of stabilizing II$_2$ by choose the edge labelled
$f$ in the defining figure as the minimal edge. Arguments as above show
that there is no problem in doing so. From henceforth, a well-chosen
edge will also be subject to this restriction.

\begin{lem}\label{triple.edges}
For any genus $g\geq 9$, a strictly optimal $G$ and a well-chosen $O(e)$, 
$G'$ contains at most one component which is a $K_{2,3}$. 
\end{lem}

\begin{proof}
Let $k$ be the number of $K_{2,3}$'s in a resulting $G'$ (note that these must
be connected components). The extra symmetry brought by these inside $G$ is 
that even if their valence two vertices are fixed, there is still a swapping 
of their stable vertices possible. We proceed to eliminate this and discuss 
the result: replace therefore each $K_{2,3}$ by a vertex, to which the edges 
incident to the original $K_{2,3}$ will be linked. This results precisely in a 
$2^k$ times decrease in the order of the automorphism group, from the 
elimination of the swapping mentioned above. However, at the same time the 
genus has dropped by $2k$. The graph $\bar{G}$ thus obtained is easily seen 
to be simle cubic, except if two edges incident to one of the collapsed 
$K_{2,3}$'s have a common endpoints, or are incident to another collapsed 
$K_{2,3}$. In the first case $O(e)$ must have been a disjoint union of stars, 
and then actually only one of them, so $g=4$; in the second case, all 
components of $G'$ must have been isomorphic; the transitivity of the action 
of $\Aut~G$ on $O(e)$ forces actually the two components linked by two edges 
in $O(e)$ to be linked by three edges in $O(e)$, and this is all of $G$. 
Then $g=6$. Both cases are outside our considerations, therefore $\bar{G}$ is 
simple cubic.

Let $\bar{g}$ be the genus of $\bar{G}$. By induction $\mu(\bar{g})\leq 9$, 
with equality only when $\bar{g}=4$; moreover $4\neq \bar{g}\leq 8$ implies 
$\mu(\bar{g})\leq 6$ and $\bar{g}\geq 9$ implies $\mu(\bar{g})\leq 3$. We are 
then studying the inequality $\mu(\bar{g})\cdot 2^{o(g-2k)+k}\leq 2^{o(g)}$ 
or equivalently $\mu(\bar{g})\leq 2^{k-l(g)+l(g-2k)}$; this is implied by 
$\mu(\bar{g})\leq 2^{o(k)}$.
\begin{itemize}
 \item $k\geq 6$ makes the last inequality strict even for $\mu(\bar{g})=9$, 
while $k\geq 4$ makes the last inequality strict for $\mu(\bar{g})\leq 6$

 \item If $\bar{g}=4$, $|\Aut~\bar{G}|\leq 72$; we need to study the cases 
$k\leq 5$. If $k=1,2$, then $g\leq 8$--too small. If $k=5$, $g=14$ and 
$|\Aut~G|\leq 72\cdot 2^5<|\Aut~C_{14}|$ so $G$ was not optimal. If $k=4$, 
$g=12$ and $|\Aut~G|\leq 72\cdot 2^4< |\Aut~C_{12}|$ so again $G$ was not 
optimal. If $k=3$, $g=10$ and since $|\Aut~G|\leq 72\cdot 2^3>|\Aut~C_{10}|$ 
we need to make use of the marking mentioned in (\ref{disc}): the three 
contracted $K_{2,3}$ must have been in the same orbit. Referring to 
a table of cubic graphs of low genus, we see that either $|\Aut~\bar{G}|=12$ 
in which case 
$|\Aut~G|<|\Aut~C_{10}|$, or $G=K_{3,3}$ so the three vertices representing 
the contracted $K_{2,3}$'s (which must not be neighbors, by the discussion on 
the simplicity of $\bar{G}$) fill one of the two sets of the partition. Then 
the marking of this partition cuts in half the number of automorphisms of 
$\bar{G}$, so $|\Aut~G|\leq 36\cdot 2^3<|\Aut~C_{10}|$ so again $G$ could 
not have been optimal.
\item If $4\neq \bar{g}\leq 8$, we need to worry about $k\leq 3$.
  \begin{itemize}
   \item  $k=1$ is only possible when $\bar{g}=7,8$. In the first situation, 
$|\Aut~G|\leq 64\cdot 2 < |\Aut~C_9|$ so $G$ was not optimal. In the second situation, the marking of the unique vertex introduced by contracting the $K_{2,3}$ is easily seen to cut at least in half the order of $\Aut~\bar{G}$; then $|\Aut~G|\leq 168\cdot 2 < |\Aut~C_{10}|$ so again $G$ was not optimal.
   \item $k=2$ is possible for $5\leq \bar{g}\leq 8$. In all cases except when $\bar{G}$ is the Petersen's graph ($\bar{g}=6$) one easily reaches the conclusion (using (\ref{low.genus.table})) that $G$ was not optimal. Using the marking (two marked points) in case of the Petersen's graph severely cuts the order of available automorphisms, to $|\Aut~G|\leq 12\cdot 2^2< |\Aut~C_{10}|$.
  \item $k=3$ and $4\neq \bar{g}\leq 8$ is again easily shown using (\ref{low.genus.table}) to lead directly (without discussing markings) to the conclusion that $G$ was not optimal.
  \end{itemize}
\item $\bar{g}\geq 9$; then by induction $\mu(\bar{g})\leq 3$ and only the cases $k=1$ and $k=2$ do not lead to the immediate conclusion that $G$ was not optimal (otherwise $o(k)\geq 2$).

When $k=2$, the graph $G'$ has at most three components, and removing one of the two $K_{2,3}$'s and its incident edges the graph is still connected, simple and cubic (since there is at most one double edge that could occur in its stabilization, by the previous reductions); the surgery would drop the genus by four, while overall $|\Aut~G|$ dropped fourfold (twice from the interchanging of the two $K_{2,3}$'s, and twice from the swapping of the two stable points of the removed $K_{2,3}$). Since $o(g)\geq o(g-4)+3$. Only now reducing the final $K_{2,3}$ would drop the genus by another two, while losing only a factor of two. Overall, we get a drop of six in genus, and a drop of eight in the order of the automorphism group. If $\bar{g}-2\geq 9$, then $o(g)\geq o(g-6)+5$ so $2^{o(g)}> \mu(\bar{G'})\cdot 2^{o(g-6)}$, so $G$ was not optimal; otherwise, we discuss as above to reach the same conclusion.
\end{itemize}
\end{proof}

To summarize the reductions so far: 

\begin{prop}\label{reduction} 
In a strictly optimal $G$ of genus $g\geq 9$, a minimal $O(e)$ may be chosen 
in such a way that $G'$ may be stabilized to a simple graph, with the 
exception of the following situations:
\begin{enumerate}
 \item $k=|O(e)|=1$ and $e$ is in the middle of a II$_+$; it is clear that 
the edges leaving from the stable ends of the II$_+$ can be at most swapped 
by $\Aut~G$, but cannot move anywhere else inside $G$. Moreover, in this case 
$G'$ is clearly connected, but stabilizing it would produce a double edge; 
furthermore, the vertices at distance one from the endpoints of the II$_+$ 
are not connected by an edge.

 \item (for $g\geq 10$)  $k=|O(e)|=1$ and this edge is incident to precisely 
one II$_2$; $G'$ is disconnected, but we know precisely an isomorphism class 
of components of $G'$ (it is easy to see that not all connected components of
$G'$ could be isomorphic; one would get two components and the genus would be 
seven, too low for our considerations). 

 \item $G'$ contains cycles; due to the length of these cycles being at least three, $|O(e)|\geq 3$, so by the discussion above, the other components of $G'$ must stabilize properly (or be cycles themselves).
 \item $G'$ may contain a unique $K_{2,3}$ as a component.
\end{enumerate}
\end{prop}

\

\noindent
As an immediate consequence, we may choose $O(e)$ minimal in such a way that 
one of the following is true:
\begin{itemize}
 \item The components of $G'$ of genus greater than two stabilize to simple
cubic graphs.

 \item $G'$ has components that are cycles, but the components that are not
stabilize to simple cubic graphs.

 \item $G'$ is connected, but stabilizing it leads to a (unique) double edge 
(this is when a unique II$_+$ occurs).

 \item $G'$ is disconnected, and one component is a II$_2$ which stabilizes
to a graph with a double edge.

 \item $G'$ contains a unique $K_{2,3}$ as a component.
\end{itemize}

\begin{rmk}(Enforcing strictness)\label{enforcing.strictness}
In certain situations, when the geometrical situation will allow, we will show that some graphs cannot be strictly optimal by the following constructions.
\begin{enumerate}
\item $G'$ is made up of three components that stabilize to simple cubic 
graphs and $O(e)$ is a disjoint union of two stars, each incident exactly 
once to each component. The two points of contact of each component with the
stars must be in the same orbit; for optimality of $G$ it is necessary that 
fixing one of the points will fix the other (i.e. the two pinched edges where 
the incidence occurs must always move together under the action of $\Aut~G$). 
Then we detach the stars, link their ends incident to each component together
to form a $K_{2,3}$, and link each of these new vertices with an edge to any 
one of the previous incidence points; we then stabilize (removing the 
pinch points at the other points of detachment). Then the automorphism group 
of the new graph is at least as large as that of the initial one, but clearly 
the minimal orbit has decreased in order; thus $G$ was not strictly optimal. 

\item $G'$ is made up of two components that stabilize to simple cubic
graphs, and $O(e)$ is a disjoint union of two or four edges. Regardless of
whether the two components of $G'$ are isomorphic, the incidence points of 
the edges in $O(e)$ with each component are in the same orbit, and their set 
is preserved by $\Aut~G$; moreover, the optimality of $G$ will dictate that 
fixing one incidence point will necessarily fix the others on its component; 
this implies that once an edge in $O(e)$ is fixed, the others will be fixed 
as well. 
  \begin{enumerate}
   \item If $k=|O(e)|=2$, then we detach the ends from one component, join 
them, and link the resulting vertex with an edge to any either of the initial 
incidence points, while stabilizing the other. 
  \item If $k=4$ we detach all the edges in $O(e)$, link the unstable vertices
of a $K_{3,3}$ with an edge removed to the endpoints of one of the edges
in $O(e)$ and stabilize the remaining endpoints.
  \end{enumerate}
  In both situations the new graph has the same genus and at least as many 
automorphisms, but the minimal orbit has strictly smaller order; thus again 
$G$ was not strictly optimal.

\item $G'$ is made up of three components, (all stabilizing to simple cubic
graphs) two of which are isomorphic, linked each by two edges to the third 
one. Then as above we may detach the edges from the two isomorphic components, 
join their free ends, and link the resulting vertices to one of the initial 
incidence vertices (stabilizing the other). As above, this is easily seen to 
contradict the strict optimality of $G$.
\end{enumerate}
\end{rmk}

\section{Proof of the Main Theorem}

We will repeatedly use an {\em exhausting subgraphs} argument. This entails 
choosing a connected component (star or edge) in $O(e)$, fixing its 
orientation (when the endpoints are in the same orbit) and then gradually 
enlarging the subgraph gotten at a certain stage by choosing one of its 
tails and adding whole components either of $O(e)$ or of $G'$ reached by that 
tail. When a component of $G'$ will be added, we will include in the new 
subgraph only the edges of $O(e)$ incident to it, and of these, in case 
$O(e)$ is a union of stars, only those that do not lead to stars whose 
center is already a vertex of the previous subgraph (in order to avoid 
cutting unnecessarily the number of tails).

At each step we look at the relative gain in the automorphism group. If a 
star is included at that step, then one of its edges is already fixed by the 
initial subgraph, then there could be at most a twofold increase in the order 
of the automorphism group at such a stage; moreover, such an increase occurs 
only when none of the vertices of the star was part of the subgraph at the 
beginning of the stage.

If however, a component is included at a certain step, then one of its 
vertices (which has valence two in $G'$) is already fixed, and that limits 
its symmetry; in other words, the automorphisms of the new subgraph fixing 
the previous one are precisely those fixing the incidence point. Once all 
these automorphisms are taken into account, all the edges incident to that 
component do not have extra freedom (they move where their incidence point 
moves), so may be added without further increase in the order of the 
automorphism group of the subgraph.

Unless otherwise noted, we will always expand the subgraphs by including 
whole components of $G'$ if the possibility exists (i.e. when not all tails 
of the subgraph gotten so far are centers of stars in $O(e)$).

During the course of the proofs, it will sometimes be convenient to disallow
certain automorphisms of a graph. In particular, sometimes we will collapse
a cycle to a vertex, but we only want to remember the automorphisms of the
resulting structure which come from the cycle. We will call the resulting 
vertex a {\em vertex with dihedral symmetry} to indicate that we do not
allow the more general automorphisms in the contracted graph.

At this time, we introduce the
following theorem, which will also be proved and used inductively in the 
course of this section:

\begin{thmb}
Suppose $g\geq 9$. If $|\Aut~C_g|>2^{o(g)}$ or $g$ is a power of two,
there is a unique strictly optimal graph of genus $g$, unless $g=10$. 
Moreover, there is a 
unique graph $G$ for which $\pi(G)=1$ in the cases $g=2^m, 3\cdot 2^m$ and
$3(2^m+2^p)$.
\end{thmb}

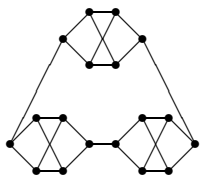
\begin{figure}[ht]
\begin{center}
\begin{picture}(80,70)
\put(5,15){\dt}
\put(5,15){\line(1,-1){10}}
\put(5,15){\line(1,1){10}}
\put(15,5){\dt}
\put(15,25){\dt}
\put(15,5){\line(1,0){10}}
\put(15,25){\line(1,0){10}}
\put(15,5){\line(1,2){10}}
\put(15,25){\line(1,-2){10}}
\put(25,5){\dt}
\put(25,25){\dt}
\put(25,5){\line(1,1){10}}
\put(25,25){\line(1,-1){10}}
\put(35,15){\dt}

\put(35,15){\line(1,0){10}}

\put(45,15){\dt}
\put(45,15){\line(1,-1){10}}
\put(45,15){\line(1,1){10}}
\put(55,5){\dt}
\put(55,25){\dt}
\put(55,5){\line(1,0){10}}
\put(55,25){\line(1,0){10}}
\put(55,5){\line(1,2){10}}
\put(55,25){\line(1,-2){10}}
\put(65,5){\dt}
\put(65,25){\dt}
\put(65,5){\line(1,1){10}}
\put(65,25){\line(1,-1){10}}
\put(75,15){\dt}

\put(25,55){\dt}
\put(25,55){\line(1,-1){10}}
\put(25,55){\line(1,1){10}}
\put(35,45){\dt}
\put(35,65){\dt}
\put(35,45){\line(1,0){10}}
\put(35,65){\line(1,0){10}}
\put(35,45){\line(1,2){10}}
\put(35,65){\line(1,-2){10}}
\put(45,45){\dt}
\put(45,65){\dt}
\put(45,45){\line(1,1){10}}
\put(45,65){\line(1,-1){10}}
\put(55,55){\dt}

\put(5,15){\line(1,2){20}}
\put(75,15){\line(-1,2){20}}
\end{picture}
\end{center}
\caption{An optimal non-$C_{10}$}
\label{genus10}
\end{figure}

\begin{rmk}
The case $g=10$ is a real exception. A graph different from $C_{10}$ with
the same number of automorphisms is depicted in Figure \ref{genus10}. In the
cases where $|\Aut~C_g|=2^{o(g)}$, non-uniqueness is the norm: a simple 
example occurs for $g=340$ ($101010100$ in binary). The candidate graph has
four ``tails'' which can be arranged in three non-isomorphic ways around
the edges of a binary tree with four ends.
\end{rmk}

We use subscripts on B in the same way as they are used on A. 

Throughout this section, we assume that Theorems A$_h$ and B$_h$ are true for 
all $h$ less than the $g\geq 9$ under consideration. 

As it relies on estimates based on the arithmetic of $g$, the proof of
the Main Theorem is somewhat tedious. Here is the outline, with details
filled in by the lemmas that occupy the rest of this section.

\begin{proof}[Proof of Theorems A and B]
As usual, the proof is divided into two cases: when $O(e)$ is a disjoint
union of stars, and when $O(e)$ is a disjoint union of edges.

{\bf Case 1:} Suppose $O(e)$ is a disjoint union of stars. Then Lemma
\ref{stars.and.genus.more.than.3} will show that if $O(e)$ consists of more 
than one star, $G$ 
is not strictly optimal. Then by the reductions of the previous section, we may
remove the star, disconnecting the graph into subgraphs of lower genus. Lemma
\ref{stars.and.genus.more.than.3} will show furthermore that the genus of 
these subgraphs is 
quite restricted. The order of the automorphism group of $G$ in this case 
will be six (for the star) times the automorphism groups of the components, so 
the stabilizations of the components must be optimal. We may then proceed by
induction: the components are of smaller genus and  we know the optimal 
pinched graphs in these genera, hence we get a bound for the automorphism
group. The second part of Theorem A about the pinched graphs follows 
similarly. Theorem B will also follow by applying it inductively to the 
components when necessary. Thus, essentially, Case I is covered by the Lemma
and induction hypotheses.

{\bf Case 2:} If $O(e)$ is a disjoint union of edges, we are in a much more
restricted situation. First of all, suppose that removing $O(e)$ and 
stabilizing results in a non-simple $G'$. We have classified the possible
behaviors in the previous section: if $G'$ remains connected, it has a unique
subgraph II$_+$. Replacing this entire subgraph by an edge, we reduce the
genus, are able to apply induction to the resulting graph (which no longer
leads to an unstable $G'$) and close this case.

If $G'$ is disconnected, we either have the case of a unique II$_2$ (touching
the minimal orbit; clearly there can be many II$_2$ in an optimal graph),
in which case we may remove $e$ and concentrate on the component which 
stabilizes well, again proving the theorems by induction, or the case of
a unique $K_{2,3}$ component in $G'$ (before stabilization). Again, we use
induction on the components of $G'$ that are not $K_{2,3}$ and obtain the
result.

Therefore, we may assume that stabilizing $G'$ does not introduce loops or
multiple edges. This finally breaks into two subcases: either $G'$ contains
cycles or it does not. Lemma \ref{all.cycles} shows that $G'$ is not a disjoint
union of cycles,
and Lemma \ref{some.cycles} then shows that a cycle in $G'$ must be unique. 
After this
is established, we see that $G$ is a collection of isomorphic subgraphs 
that stabilize well arranged around a cycle. Again we may apply induction to 
the subgraphs to obtain the theorems.

If $G'$ does not contain a cycle, then the orbit of a minimal edge must be
small, and the estimates finish the work. The details are recorded in
\ref{isolated.edges.and.genus.more.than.3}.

The proofs of the lemmas will show that if $g=2^m$ there is only one optimal 
graph $G$ for which $\pi(G)=1$. If $g=3\cdot 2^m$, then $A_m$ is optimal, but
$\pi(A_m)=\frac12$. This part of Theorem B will be show by showing that the
 only other graphs with at least $2^{o(g)}$ 
automorphisms in these genera are those linking a $B_{m+1}$ and a $B_m$ by an 
edge;
these graphs are not optimal (they have exactly $2^{o(g)}$ automorphisms), 
but they satisfy $\pi=1$. A similar observation
is true in the case $g=3(2^m+2^p)$. This proves the last part of Theorem B,
which is essential in the induction.
\end{proof}

\subsection{$G'$ has components that are cycles}

\begin{lem}\label{all.cycles}
If $G'$ is a disjoint union of cycles, then $G$ cannot be optimal.
\end{lem}
\begin{proof}
{\bf Case 1}:
If $G'$ is a single cycle, and $k=|O(e)|$, then the length of the cycle 
is $2k$ with $4k$ automorphisms; $g=k+1$ and $2^{k+1}>4k$ for $k\geq 5$, 
i.e. for $g\geq 6$, which is certainly the case; such a $G$ could not be 
optimal. 

{\bf Case 2:} If $G'$ is disconnected and $O(e)$ is a disjoint union of $k$ 
stars, then all the connected components of $G'$ are isomorphic (since all 
the vertices of the stars are in the same orbit). Then $G'$ having a 
connected component which is a cycle means that $G'$ is a disjoint union of 
$s$ cycles of the same length, say $t$; note that $t\geq 3$, otherwise $G$ 
could not have been simple. We note that the transitivity of the action of 
$\Aut~G$ on the orbit $O(e)$, coupled with the disconnectedness of $G'$, 
prohibits the stars from having more than one contact point with any 
component of $G'$.  The genus of $G$ is $g=2k+1$, and we have the ``contact 
formula'' $st=3k$. In estimating $|\Aut~G|$ we may simply replace the cycles 
in $G'$ by (contract them to) vertices of valence $t$ with dihedral 
symmetry. We get a graph with $k$ vertices of valence three and $s$ vertices
of valence $t$.

We use growing trees: choose a star ($k$ choices), fix its edges (at most six 
choices), then expand this tree and subsequent trees by reaching, from a 
tail, to adjacent vertices not included in the tree so far. Due to the at most 
dihedral symmetry around each vertex of the contracted graph, at each step 
in the process the size of the automorphism group increases at most two-fold; 
this occurs precisely when a tail of order three in $G$ has only one neighbor 
in the tree gotten so far (say this happens $a$ times), or when a tail of 
order $t$ in $G$ has at most two neighbors (one, if $t=3$) in the tree (say 
this happens $b$ times). Thus $|\Aut~G|\leq 6k\cdot 2^{a+b}$ and (denoting 
by $n$ the number of vertices of the contracted graph) $n=s+k\geq 4+2a+(t-2)b$ 
(if $t\geq 4$) respectively $n\geq 4+2a+2b$ (if $t=3$). In the first instance 
we get $a\leq \frac{s+k}{2}-2-\frac{t-2}{2}b$, so 
$|\Aut~G|\leq \frac{3}{2}k\cdot 2^{\frac{s+k}{2}-b\frac{t-4}{2}}$ (and note 
that $s\leq \frac{3k}{4}$); in the second instance we get directly 
$a+b\leq \frac{s+k}{2}-2$ (but $s=k$) so $|\Aut~G|\leq \frac{3}{2}k\cdot 2^k$. 
Both times we compare to $2^{o(2k+1)}=2^{2k+1-l(2k+1)}\geq2^{k+o(k)}$ 
(since $l(2k+1)\leq l(k)+1$; see (\ref{estimate.l.function})).

Now both inequalities are implied by $\frac{3}{2}k\leq 2^k$ which is strict 
by (\ref{daily.inequalities}), since $k\geq o(k)+1\geq o(k+1)$.

{\bf Case 3}: If $G'$ is disconnected and $O(e)$ is a disjoint union of 
isolated edges, then there are at most two isomorphism classes of connected 
components (according to whether both endpoints of $e$ are in the same orbit 
or not). 

We will deal directly with the general case, in which there are two 
isomorphism classes of cycles in $G'$; the simpler case may be dealt with by 
taking $s_1=s_2$ below; all the estimates still work then. Thus there are 
$s_1$ cycles $H_1$ incident by $t$ edges to each of $n_1$ neighbors (these 
have length $n_1t$), and $s_2$ cycles $H_2$ incident by $t$ edges to each of 
$n_2$ neighbors (these have length $n_2t$). Then $k=s_in_it$ and $g=k+1$.

We need to show first that any two cycles may be linked by at most an edge in 
$O(e)$, i.e. $t=1$.

If two cycles would have $t\geq 3$ common edges in $O(e)$ incident to both, 
then clearly fixing one cycle will fix the whole graph. Then 
$|\Aut~G|\leq 2s_1n_1=2k<2^{o(k+1)}$ for $k\geq 5$ using 
(\ref{daily.inequalities}). But $k\leq 4$ and $t\geq 3$ is only possible when 
$G$ is made up of two isomorphic cycles sharing an edge between their 
vertices (one from each component); in that case the genus is too small (five 
or six) for our considerations. Thus such a $G$ could not be optimal.

If two cycles would have precisely two edges in common, then 
$|\Aut~G|\leq 2k\cdot 2^{a+b}$, where $a$ is the number of times we get an 
involution of a cycle of type $H_1$ by including it when expanding the 
subgraph at a tail incident to it; only one neighbor of this $H_1$ must have 
been included in the expanding subgraph previously, so $n_1-1$ new cycles of 
type $H_2$ will be incident to the newly increased subgraph afterwards; and 
similarly for $b$. Thus $s_2\geq n_1+a(n_1-1)$ and $s_1\geq 1+b(n_2-1)$. 

Note that if $n_1=1$ we get $s_2=1$; then $t=2$ would force the existence of 
double edges in $G$. Thus $\min (n_1,n_2)\geq 2$; as a consequence, since 
$t=2$, $s_i\leq \frac{k}{2}$.

If $n_1\geq 3$, then 
$a\leq \frac{s_2-1}{n_1-1}-1\leq \frac{s_2-1}{2}-1< \frac{k}{4}-1$; similarly 
$b\leq \frac{s_1-1}{n_2-1}<s_1=\frac{k}{n_1}\leq \frac{k}{3}$. Then overall 
$|\Aut~G|\leq 2k\cdot 2^{\frac{7k}{12}-1}<k\cdot 2^{\frac{k}{2}}$; now, since 
$k=2s_1n_1\geq 12$, this is easily seen to be strictly less than $2^{o(k+1)}$ 
by (\ref{daily.inequalities}) so again $G$ is not optimal.

If $n_1=n_2=2$, then $G$ is formed of cycles of length four, each with two 
neighbors with which it is linked by two edges; it is immediate that the 
edges between two adjacent cycles should be linked at their opposite vertices 
for maximum symmetry gain. These cycles are actually isomorphic, and there 
are precisely $\frac{k}{2}$ of them in $G$. However, in this case an 
involution in one cycle will force an involution in a neighboring cycle; 
overall $|\Aut~G|\leq k\cdot 2^{\frac{k}{2}}<2^{o(k+1)}$ as soon as 
$k\geq 10$, so $G$ is not optimal.

We note that for $g=9$, i.e. $k=8$, 
$|\Aut~G|\leq 8\cdot 2^4=2^{o(g)}<|\Aut~C_9|$ so $G$ cannot be optimal.

We also note that $g=2^m$ is not possible in the above configuration.

Thus we have reduced to only the possibility $t=1$.

We can estimate as above by first collapsing all cycles to vertices, then by 
expanding trees. So we get a graph $\bar{G}$ with $s_1$ vertices of order 
$t_1$ and $s_2$ vertices of order $t_2$, with dihedral symmetry around each 
vertex, and, most importantly, with $|\Aut~\bar{G}|=|\Aut~G|$. Note that 
$k=|O(e)|=s_1t_1=s_2t_2$, $t_i\geq 3$ and $g=k+1$. Without loss of generality,
 we may assume $t_1\geq t_2$.

Construct the exhausting trees by first choosing an edge in $O(e)$ ($k$ 
choices), then fixing its orientation (at most two choices); 
afterwards, each tail which has no more than two neighbors in the trees 
constructed so far may bring at most an extra involution among the edges 
ending at it (since two of those are fixed necessarily); in case $t_i=3$, 
there should be only one neighbor of that tail among the vertices touched by 
the tree so far in order for that involution to exist. Say this situation 
occurs $a$ times for the cycles of length $t_1$ and $b$ times for the cycles 
of length $t_2$. Overall $|\Aut~G|\leq 2k\cdot 2^{a+b}$. Then we have the 
following estimates:

\begin{enumerate}
\item $s_1+s_2\geq 2+a(t_1-2)+b(t_2-2)$ when $t_i\geq 4$; thus $b\leq \frac{s_1+s_2-2}{t_2-2}-a\frac{t_1-2}{t_2-2}$, so $|\Aut(\bar{G})|\leq 2k\cdot 2^{\frac{s_1+s_2-2}{t_2-2}}$
\item $s_1+s_2\geq 2+a(t_1-2)+2b$ when $t_1>t_2=3$; thus $b\leq \frac{s_1+s_2-2}{2}-a\frac{t_1-2}{2}$, so $|\Aut(\bar{G})|\leq k\cdot 2^\frac{s_1+s_2}{2}$. 
\item $s_1+s_2\geq 2+2a+2b$ when $t_1=t_2=3$ (so $s_1=s_2$), so $|\Aut(\bar{G})|\leq k\cdot 2^{s_1}$.
\end{enumerate}

In all instances we compare with $2^{o(k+1)}$.

In the first subcase, $s_1=\frac{k}{t_1}\leq \frac{k}{t_2}=s_2$, so the inequality to prove is implied by $2k\cdot 2^\frac{2(k-t_2)}{t_2(t_2-2)}\leq 2^{o(k+1)}=2^{k+1-l(k+1)}$; this in turn is implied, using (\ref{estimate.l.function}) and the fact that $\frac{k-t_2}{t_2(t_2-2)}\leq \frac{k-4}{8}$, by $4k\sqrt{k+1}\cdot 2^\frac{k-4}{4}\leq 2^{k+1}$. This last one is equivalent to $k\sqrt{k+1}\leq 2^\frac{3k}{4}$, which is easily seen to be strict for $k\geq 5$; however, $g\geq 8$ implies $k\geq 7$ so in this first subcase we always get a strict inequality.

In the second subcase, $s_1=\frac{k}{t_1}\leq \frac{k}{3}=s_2$; as above, we reduce to $k\sqrt{k+1}\leq 2^\frac{2k}{3}$, which is easy to be seen as true (strict inequality) for $k\geq 6$, which again is what we needed given that $k\geq 8$. 

The third subcase is reduced to the same inequality as the second, so we are done.
\end{proof}

\begin{lem}\label{some.cycles}
If $G$ is a strictly optimal graph, $G'$ cannot contain more than one 
component which is a cycle. Moreover, Theorem B holds for those $g$ for which
$C_g$ contains an isolated cycle.
\end{lem}

\begin{proof}
Due to (\ref{all.cycles}), we are left to analyze the case where some
components of $G'$ are cycles while the others are not.

First, (\ref{reduction}) allows us to assume that the $s_2$ components 
$H_i$ of $G'$ which are not cycles may be stabilized without problems (no 
double edges or loops occur); in particular, the arithmetic genus, denoted 
by $h$, of these components is at least three. Then the induction hypothesis 
shows that $|\pi(H_i)|\leq 2^{o(h)}$ for $h\geq 3$. $G'$ has also $s_1$ 
components which are cycles, each with $n_1$ neighbors (components of $G'$ at
distance one); we denote by $t$ the {\em incidence degree} of two components 
in different isomorphism classes, i.e. the number of edges in $O(e)$ joining 
two such components; we also denote by $n_2$ the number of neighbors of a 
component $H_i$. We have: the length of the cycles is $n_1t$, 
$k=s_1n_1t=s_2n_2t$ and $g=s_2(h-1)+1+k$. Note that $n_1t\geq 3$ (to avoid 
double edges in $G$) and since we want $s_1=1$ it is enough to show that 
$k\geq 6$ leads to contradictions. We will discuss also what happens when 
$s_1=1$.

Estimate $|\Aut~G|$ by first contracting the cycles to vertices with dihedral 
symmetry at the edges around them, then using again an exhausting subgraphs 
argument; we start by choosing a cycle ($s_1$ choices), then fixing its 
orientation ($2n_1t$ choices). We get 
$|\Aut~G|\leq 2k\cdot  2^{a}\cdot 2^{s_2o(h)}$, and would like to compare 
this to $2^{o(g)}$. Due to the way we mentioned we expand the 
subgraphs, preferably at their tails incident to components $H_i$ (when 
they exist at a certain stage) we see that the only possibility of gaining 
extra symmetry when forced to incorporate a cycle is when that cycle had at 
most two incident edges (and precisely one if its length $n_1t=3$). Including 
such a cycle will immediately yield $n_1t-2$ (respectively $2$) tails 
incident to components $H_i$; since $t$ was the incidence, we see that we 
have: if $t\geq 2$, only one component $H_i$ was incident to this cycle, so 
$n_1-1$ new components will be reached by tails after the cycles is included 
in the newly increased subgraph; if $t=1$ but $n_1\geq 4$, at least $n_1-2$ 
new components will be reached; and finally if $t=1$ but $n_1=3$, exactly 
$2=n_1-1$ new components will be reached. This happens each of the $a$ times. 
Thus we have:

\begin{itemize}
 \item $s_2\geq n_1+a(n_1-1)$ if either $t\geq 2$ or $t=1,n_1=3$
 \item $s_2\geq n_1+a(n_1-2)$ if $t=1$ and $n_1\geq 4$ 
\end{itemize} 

 Quick manipulations lead us to the inequality 

$$2k\leq 2^{s_2h-s_2+k+1-l(s_2h-s_2+k+1)-a-s_2h+s_2l(h)}$$

Using (\ref{estimate.l.function}) we get 
$l(s_2h-s_2+k+1)\leq l(s_2)l(h)+l(k+1-s_2)$ so the above inequality is 
implied by $2k\leq 2^{B(s_2,h)-a+o(k+1-s_2)}$. 

If $s_2=1$, then $n_1t\geq 3$ (the cycles must have length at least three, 
otherwise $G$ would have double edges); moreover, the cycles are connected to 
the ``core'' component $H$ of $G'$ by edges starting from all of their 
vertices. Thus (dispensing with the above bound on $|\Aut~G|$) we have in 
fact $|\Aut~G|=|\Aut~H|$. But $g=h+k=h+s_1n_1t\geq h+3$, so 
$o(g)\geq o(h)+o(3)=o(h)+2$; now either $h\geq 9$, in which case the 
induction hypothesis says that 
$|\Aut~H|\leq 3\cdot 2^{o(h)}<4\cdot 2^{o(h)}\leq 2^{o(g)}$ (so $G$ could not 
have been optimal), or $h\leq 8$, in which case the table 
(\ref{low.genus.table}) and the estimates (\ref{growth.candidates}) show 
again that $G$ could not have been optimal. Thus this case cannot occur for 
an optimal graph $G$.

Thus from now on $s_2\geq 2$, which implies $n_1\geq 2$. 

Suppose $n_2\geq 2$ (or, equivalently, $s_1\geq 2$).

If $n_1\geq 4$ (and any $t$), then $s_2\geq 4$ (so by (\ref{main.inequality})  
$B(s_2,h)\geq \lfloor \frac{s_2+1}{2}\rfloor$) and 
$s_2\leq \frac{k}{n_2t}\leq \frac{k}{2}$;  also 
$a\leq \lfloor \frac{s_2-1}{2}\rfloor-1$; then $B(s_2,h)-a\geq 2$. The 
inequality is then implied by $2k\leq 2^{2+o(\lceil \frac{k}{2}\rceil +1)}$, 
or equivalently $\frac{k}{2}\leq 2^{o(\lceil \frac{k}{2}\rceil +1)}$. Now 
(\ref{daily.inequalities}) shows that this is strict, so $G$ could not have 
been optimal.

If $n_1=3$ (and any $t$), then $s_2\geq 3$, $k=3s_1\geq 6$ and $a\leq \lfloor \frac{s_2-1}{2}\rfloor-1$; in the same time $s_2=\frac{k}{n_1t}\leq \frac{k}{3}$. (\ref{main.inequality}) (more precisely, (\ref{long.table})) shows that $B(s_2,h)-\lfloor \frac{s_2+1}{2}\rfloor\geq -1$ (with equality if and only if $s_2=3$ and $l(h)=1$) so $B(s_2,h)-a\geq 1$; the inequality is then implied by $2k\leq 2^{1+o(\lceil \frac{2k}{3}\rceil +1)}$ which is strict by (\ref{daily.inequalities}) so again $G$ is not optimal.

Thus we must have both $s_1=n_2=1$ and $s_2=n_1\geq 2$ in an optimal $G$; the 
lemma is proved.

Then $a=0$ and $k=s_2t$; we must have $k\geq 3$ since that is the length of 
the cycle.

If $t\geq 2$, $2\leq s_2\leq \lfloor \frac{k}{2}\rfloor$ and $B(s_2,h)\geq l(h)\geq 1$ (using (\ref{long.table}) in (\ref{main.inequality})). The inequality is implied by $k\geq 2^{o(\lceil \frac{k}{2}\rceil+1)}$, which by (\ref{daily.inequalities}) is strict except for $k=2,4,8$. But $k=s_2t\geq 4$ so only $k=4,8$ need consideration. $k=4$ may occur only when $s_2=t=2$ (and $l(h)=1$) and then (\ref{enforcing.strictness}) shows that $G$ was not optimal. $k=8$ may occur for either $s_2=2,t=4$ or $s_2=4,t=2$; however the inequality $2k\leq 2^{B(s_2,h)+o(k+1-s_2)}$ is easily seen to be strict in these cases, so again $G$ was not optimal.

Then we must have $t=1$ so $k=s_2\geq 3$; the inequality becomes $2k\leq 2^{B(k,h)}$. Then (\ref{main.inequality}) and (\ref{daily.inequalities}) show that the inequality is strict for all $k\geq 6$.

If $k=5$, the inequality becomes $10\leq 2^{B(5,h)}$ which is strict for $l(h)\geq 2$ by (\ref{long.table}), and false for $l(h)=1$

If $k=4$, the inequality becomes $8\leq 2^{B(4,h)}$ which is again strict for $l(h)\geq 2$ by (\ref{long.table}), with equality for $l(h)=1$.

If $k=3$, the inequality becomes $6\leq 2^{B(3,h)}$ which is strict for 
$l(h)\geq 3$ by (\ref{long.table}), and fails for $l(h)=1$ or when $l(h)=2$ 
and $l(3h)=4$. 

Even if $g$ would be $9$, $2^n$, or $3\cdot 2^n$, a small calculation shows
that the graphs with a cycle do not give an optimal graph.

The cases which remain are exactly the $g$ for which $C_g$ contains an 
isolated cycle, so Theorem B holds in these cases by induction.
\end{proof}

From now on, we may assume that $G'$ contains no cycles. We first address
the case that the minimal orbit is a disjoint union of stars.

\begin{lem}\label{stars.and.genus.more.than.3}
Let $G$ be a simple cubic graph of genus $g\geq 9$, with a minimal orbit 
$O(e)$ a disjoint union of $k$ stars, with all the components of $G'$ of 
genus $h\geq 3$ (by induction on Theorem A, $\pi(G_i)\leq 2^{o(h)}$). Then 
$G$ is not strictly optimal as soon as $k\geq 2$; moreover, for $k=1$, $G$ is 
strictly optimal only if either $l(h)=1$ or $h=3\cdot (2^m+2^p)$ with 
$|m-p|\geq 5$. Therefore, by induction, Theorems A holds in these cases.
\end{lem}

\begin{proof}
The reduction (\ref{reduction}) shows that in case $O(e)$ is a disjoint union of stars, the components of $G'$ must stabilize without problems; thus we may use the induction hypothesis in this case.

If $G'$ is connected, then 
$|\Aut~G|\leq |\Aut~G^\mrm{'stab}|\leq |\Aut~C_{g'}|$ where $g'=g-2k$. Using 
(\ref{growth.candidates}) we see that $G$ could not have been optimal for any 
$k\geq 1$.

Then let $s\geq 3$  be the number of components of $G'$ (if $s\leq 2$ then 
$G'$ is connected, as a star cannot be incident to a given component twice 
without having actually all tails in that component; thus $G'$ disconnected 
implies that each star is incident to three distinct component of $G'$). Let 
$t$ be the number of edges in $O(e)$ incident to a given component. Then 
$3k=st$, and $g=sh-s+2k+1$. 

We also note that $t=1$ implies $k=1$ (the subgraph made up of a star and the three components to which it is incident would be a connected component of $G$, thus the whole $G$). 

Using the exhausting subgraphs argument we get $|\Aut~G|\leq 6k\cdot 2^{s\cdot o(h)}\cdot 2^a$, where $a$ is the number of times we might have gained a twofold increase in the automorphism group of the subgraph by including a (new) star at a tail incident to it. Due to the way we construct these enhausting subgraphs, each inclusion of a star counted among the $a$ ones will make the new subgraph incident to two other components to which the previous subgraph was not incident. Thus we see that $s\geq 3+2a$ (since at the beginning we already had a star incident to three components).

Thus we would like to show that $6k\cdot 2^{s\cdot o(h)+\integer{\frac{s-3}{2}}}\leq 2^{o(g)}$, or equivalently (using (\ref{estimate.l.function}))  $$3k\leq 2^{o(2k-s+1)+sl(h)-l(sh)-\integer{\frac{s-1}{2}}}=2^{A(s,h)+o(2k-s+1)}$$

If $t\geq 3$ then $s\leq k$ so the inequality is implied by $3k\leq 2^{A(s,h)+o(k+1)}$; since $s\geq 3$, by (\ref{main.inequality}) $A(s,h)\geq 1$ and using (\ref{daily.inequalities}) this is easily seen to be strict for all $k\geq 1$, so such a $G$ cannot be optimal.

If $t=2$, then setting $k=2u$ we get $s=3u$; we need to study when 
$6u\leq 2^{o(u+1)+A(3u,h)}$. For $u\geq 2$, $A(3u,h)\geq 2$ by (\ref{main.inequality}) and then again (\ref{daily.inequalities}) shows that the inequality is strict; such $G$ cannot be optimal. If $u=1$, there are three connected components in $G'$, isomorphic and linked by two stars; it is apparent that already 
$M(G)\geq 3$. The inequality becomes $6\leq 2^{o(2)+3l(h)-l(3h)-1}=2^{3l(h)-l(3h)}\geq 2^{l(h)}$. Then clearly for $l(h)\geq 3$ the inequality is strict, as it is for $l(h)=2$ but $l(3h)\leq 3$. Thus either $l(h)=1$ or  $l(h)=2$ and $l(3h)=4$; however, even in these cases (\ref{enforcing.strictness}) shows that $G$ could not have been strictly optimal. 

We are left with considering the case $t=1$. Then $k=1$ as remarked before, 
so the inequality becomes $6\leq 2^{3l(h)-l(3h)}$. From (\ref{long.table}), 
we see that the inequality is strict as 
soon as $l(h)\geq 3$. Moreover, for $l(h)=2$ and $l(3h)\leq 3$ we get again 
a strict inequality. 

Thus only $l(h)=1$ or $l(h)=2,l(3h)=4$ are left. In all cases, $M(G)=3$, and overall $|\Aut~G|\leq 3\cdot 2^{o(g)}$ as claimed; moreover, equality may occur only when $h=3\cdot 2^m$. We only need to show that any graph with more than $2^{o(g)}$ automorphisms is forced to have $M(G)\geq 3$, and it is clear that the above reductions prove just that, so we are done.
\end{proof}

Finally, the case of $O(e)$ a disjoint union of edges must be analyzed.

\begin{lem}\label{isolated.edges.and.genus.more.than.3}
Let $G$ be a simple cubic graph with a minimal orbit $O(e)$ a disjoint union 
of $k$ edges, with all components of $G'$ of genus at least three 
and stabilizing to simple cubic graphs.
Then $G$ is not strictly optimal as soon as 
$k\geq 5$. Moreover, Theorems A and B hold in these cases.
\end{lem}

\begin{proof}
{\bf Case 1: $G'$ is connected} and (according to the hypothesis) it 
stabilizes without double edges or loops.

If $g=9$ then $g'=g-k\leq 8$; it is clear that the inequality $|\Aut~G|\leq|\Aut~G'|$ and the table (\ref{low.genus.table}) show that $G$ could not be optimal for any $k$.

If $g=8$ then $g'\leq 7$; $2^{o(g)}=128$ and it is clear that $|\Aut~G|\leq|\Aut~G'|$ and the table (\ref{low.genus.table}) show that $G$ could not be optimal for any $k$.

If $g=2^m\geq 16$, then for $k\geq 2$ we get $|\Aut~G|\leq|\Aut~G'|\leq |\Aut~C_{g'}|$ and (\ref{growth.candidates})) shows that $G$ could not be optimal for such $k$; and if $k=1$, then $l(g')\geq 2$ so by induction the Main Theorem shows that $|\Aut~C_{g'}|\leq 3\cdot2^{g'-2}=3\cdot 2^{g-3}<2^{o(g)}$; thus again $G$ could not have been optimal.

Otherwise, 
$|\Aut~G|\leq|\Aut~G'|\leq |\Aut(G^{'\mrm{stab}})|\leq |\Aut(C_{g'})|$. 
Now $g'=g-k$ so (\ref{growth.candidates}) shows that $G$ is not optimal as 
soon as $k\geq 2$. If $k=1$, (\ref{growth.candidates}) gives 
$|\Aut~G|\leq |\Aut~G^{'\mrm{stab}}|\leq |\Aut~C_{g'}|$; we would like to 
compare the last term with $2^{o(g)}\geq 2^{o(g')}$. If 
$|\Aut~G'|\leq 2^{o(g')}$, we get $|\Aut~G|\leq 2^{o(g)}$ which would prove
Theorem A in this case. If however $|\Aut~G^{'\mrm{stab}}|>2^{o(g')}$, then 
$M(G^{'\mrm{stab}})\geq 3$ by induction. Now $k=1$ shows that the edges of 
$G^{'\mrm{stab}}$ which are pinched by $e$ must be in an orbit by themselves 
under the action of $\Aut~G'$. Thus $\Aut~G'$ is the stabilizer of the two 
unstable paths (edges pinched by $e$). If $M(G')=4$ we then get at least a 
twofold reduction in $|\Aut~G^{'\mrm{stab}}|$ vs. $|\Aut~G'|$, while if 
$4\neq M(G') \geq 3$ we get at least a threefold reduction in the same 
comparison. The induction shows that $M(G')=4$ happens only for 
$g'=g-1=5\cdot 2^m+1$ (and then $|\Aut~G|<2^{o(g)}$) and otherwise 
$|\Aut~G'|\leq 2^{o(g')}$ as desired, with equality if and only if  
$g'=g-1=9\cdot 2^m+u$ ($u=0,2$); if $u=2$, $o(g)>o(g')$ so again $G$ could
not have been optimal, while if $u=0$ the induction shows that the star 
disconnecting the $G^{'\mrm{stab}}$ is the minimal orbit; inserting $e$ 
produces actually a sixfold loss of symmetry, so again $G$ could not have 
been optimal since in any case we get $|\Aut~G|< 2^{o(g)}$. Note that the 
reasoning applies when $g=9$ or $g=2^m$ with $m\geq 4$ to yield that such 
a $G$ could not have been optimal (strict or not).

{\bf Case 2: $G'$ is disconnected and its components are split in two 
isomorphism classes}, one of $s_1$ subgraphs isomorphic to $H_1$ and one with $s_2$ subgraphs isomorphic to $H_2$. Say each $H_1$ is at distance one from $n_1$ components isomorphic to $H_2$; define $n_2$ similarly. Define also $t$ to be the number of edges linking two neighboring components $H_1$ and $H_2$. Note that $g=s_1(h_1-1)+s_2(h_2-1)+k+1$, $k=s_1n_1t=s_2n_2t$; note also that $s_1=1$ if and only if $n_2=1$, and similarly $s_2=1$ if and only if $n_1=1$. We may assume, without loss of generality, that $s_1\leq s_2$.

Start again with an exhausting subgraphs argument; the initial step is choosing an edge and its orientation (if all components of $G'$ are isomorphic). We get $|\Aut~G|\leq k\cdot 2^{s_1o(h_1)+s_2o(h_2)}$ and we would like to study the inequality $k\cdot 2^{s_1o(h_1)+s_2o(h_2)}\leq 2^{o(g)}$. This reduces to $$k\leq 2^{s_1l(h_1)+s_2l(h_2)+k+1-s_1-s_2-l(s_1h_1+s_2h_2+k+1-s_1-s_2)}.$$ This is implied, via (\ref{estimate.l.function}), by the inequality 

\begin{equation*}
\tag{*}
k\leq 2^{B(s_1,h_1)+B(s_2,h_2)+o(k+1-s_1-s_2)}.
\end{equation*}

Assume first $t\geq 2$ and any $g\geq 9$.

If $n_1,n_2\geq 2$, then $s_i\geq 4$ so by (\ref{main.inequality}, 
(\ref{long.table})) we see that $B(s_i,h_i)\geq 2$; moreover, 
$s_i= \frac{k}{n_it}\leq \frac{k}{4}$ so then the inequality is implied by 
$\frac{k}{16}\leq 2^{o(\lceil \frac{k}{2}\rceil+1)}$ which is strict for all 
$k\geq 1$ in light of (\ref{daily.inequalities}).

If $n_1\geq 2$ and $n_2=1$, then actually $s_1=1$, $s_2=n_1\geq 2$ and 
$k=s_2t\geq 4$. Then $B(s_2,h_2)\geq 1$ and $s_2=\frac{k}{t}\leq \frac{k}{2}$ 
so the inequality is implied by 
$\frac{k}{2}\leq 2^{o(\lceil \frac{k}{2}\rceil+1)}$ 
which is strict for all $k\geq 1$ by (\ref{daily.inequalities}).

If finally $n_1=n_2=1$ then $s_1=s_2=1$, $k=t$ and the inequality becomes $k\leq 2^{o(k-1)}$. This is strict for $k\geq 5$ by (\ref{daily.inequalities}), so we only need to consider the cases when $G$ is made up of two (non-isomorphic) components linked by two, three or four edges in the same orbit of $\Aut~G$.

If $g=9$, then if $k=4$ the two components must have genus three each, so overall $|\Aut~G|\leq 4\cdot 2^4$ (since $\mu_1(3)=1$); if $k=3$ then one component has genus three and the other genus four, so overall $|\Aut~G|\leq 3\cdot 2^5$ (since $\mu_1(3)=\mu_1(4)=1$); and if $k=2$ we may have either two components of genus four each, or a component of genus three and one of genus five, in either case getting $|\Aut~G|\leq 2\cdot 2^5$; then (\ref{low.genus.table}) shows that $G$ could not have been optimal.

If $g=8$, then only $k\leq 3$ is possible, otherwise one of the components will have genus less than two. If $k=3$, the components are tetrahedra, so $|\Aut~G|\leq 3\cdot 2^4$; if $k=2$ then one component is a tetrahedron and the other one has genus four, so $|\Aut~G|\leq 2\cdot 2^5$; in both cases we get less than $2^7=2^{o(g)}$ so such a $G$ cannot be optimal.

If $g=2^m\geq 16$, then $|\Aut~G|\leq k\cdot 2^{o(h_1)+o(h_2)}\leq k\cdot 2^{h_1+h_2-2}=k\cdot 2^{g-k-1}<2^{g-1}=2^{o(g)}$, so $G$ cannot be optimal. 
Similar estimates show that such a $G$ is not optimal if $g=3\cdot 2^m$ or
$g=3(2^m+2^p)$ (with the usual restriction on $m$ and $p$).

Otherwise, for $g>9$ which is not a power of two, 
(\ref{enforcing.strictness}) reduces the search for strictly optimal $G$'s to 
the case $|O(e)|=3$. We note also that, since $g=h_1+h_2+2$ and 
$l(g)\leq l(h_1)+l(h_2)+1$, 
$|\Aut~G|\leq 3\cdot 2^{o(h_1)+o(h_2)}=3\cdot 2^{o(g)-(l(h_1)+l(h_2)+2-l(g))}\leq \frac{3}{2}2^{o(g)}$; 
we also have $M(G)=3$.  Fixing an edge in $O(e)$ must have as a result the
 fixing of all the edges in $O(e)$; regarding the components, fixing one of 
the three incidence points must fix all three (growing the subgraphs by 
choosing an edge, then including a component will automatically fix all the 
edges, therefore the other component's incidence points). Then, on one side, 
this implies $M(H_i^\mrm{stab})\geq 3$; on the other it 
says  $\pi(H_i)\leq \frac{1}{3}|\Aut~H_i|$. By induction we know that only 
for $h_i=9\cdot 2^{n_i}+u_i$ ($u=0,1,2$) may we get the coefficient 
$\mu(H_i)=3$; otherwise $\mu(H_i)\leq \frac{3}{2}$, so 
$\pi(H_i)\leq \frac{1}{2}2^{o(h_i)}$ (for at least one $i$) so it is easy to 
see then that $|\Aut~G|\leq \frac{3}{4}2^{o(g)}$ so $G$ would not be optimal. 
However, a quick computation shows that even if $h_i=9\cdot 2^{n_i}+u_i$ one 
still gets $l(h_1)+l(h_2)\geq l(g)$ so $|\Aut~G|<2^{o(g)}$ so $G$ could not 
be optimal in this case.

Thus we must have $t=1$ for a strictly optimal $G$. 

If $n_1,n_2\geq 3$, then $s_i\geq 3$ so again $B(s_i,h_i)\geq 2$, 
$s_i\leq \frac{k}{3}$ and the inequality is implied by 
$\frac{k}{16}\leq 2^{o(\lceil \frac{k}{3}\rceil+1)}$, which is strict for 
all $k\geq 1$ by (\ref{daily.inequalities}).

If $n_1\geq 3$ and $n_2=2$ (similarly for $2=n_1<n_2$), then 
$s_2\geq 3$, $s_1\geq 2$ so $B(s_1,h_1)\geq 1$, $B(s_2,h_2)\geq 2$, while 
$s_2\leq \frac{k}{2}$ and $s_1\leq \frac{k}{3}$. Then the inequality is 
implied by $\frac{k}{8}\leq 2^{o(\lceil \frac{k}{6}\rceil)+1}$, again strict 
for all $k\geq 1$ by (\ref{daily.inequalities}).

If $n_1=n_2=2$, then $G$ is a pseudocycle formed by components linked to each 
of their two neighbours by an edge, with the alternate components isomorphic; 
then $s_1=s_2$ (note that $g$ is odd so $g=2^m$, $3\cdot 2^m$, and 
$3(2^m+2^p)$ with $m, p>0$ are impossible, and also that 
$g=9$ would force all components to have genus two, which is again ruled out). 
If $s_1\geq 3$ the cycle has length $2s_1\geq 6$ and cannot be strictly optimal: grouping adjacent components two by two (linked by the common edge), pinching this edge and finally linking the new $s_1$ components to the vertices of a cycle of length $s_1$ produces a graph with the same number of automorphisms, but with fewer edges in the minimal orbit.
If $s_1=2$ the inequality becomes $4\leq 2^{l(h_1)+l(h_2)}$, so as soon as $\max(l(h_1),l(h_2))\geq 2$ we get a strict inequality. Thus for an optimal $G$ one must have $l(h_i)=1$. Then $g=2h_1+2h_2+1$ and $|\Aut~G|\leq 4\cdot 2^{2o(h_1)+2o(h_2)}=2^{o(g)-(2l(h_1)+2l(h_2)+1-l(2h_1+2h_2+1))}$; this is strictly less than $2^{o(g)}$ by (\ref{estimate.l.function}) so $G$ could not have been optimal.

If $n_1=1$ and $n_2\geq 3$, then $s_2=1$ and $s_1=n_2=k\geq 3$; the 
inequality (*) becomes $k\leq 2^{B(k,h_1)}$.
\begin{itemize}
 \item If $l(h_1)\geq 2$, (\ref{long.table}) in the proof of 
(\ref{main.inequality}) and (\ref{daily.inequalities}) show that 
$2^{B(k,h_1)}\geq 2^{\lfloor \frac{k+1}{2}\rfloor}\geq k$, 
with at least one of these inequalities strict (since $k\geq 3$); thus
$G$ could not have been optimal.

  \item If $l(h_1)=1$ a similar argument shows that only $k=3$ should be 
considered in all cases. If $g=9$ then $k=3$ is not possible (a component 
would have genus less than two). If $g=2^m$, $3\cdot 2^m$, or $3(2^m+2^p)$ 
and $k=3$, $3h_1+h_2=g$ so 
$|\Aut~G|\leq 3\cdot 2^{3o(h_1)+o(h_2)}= 3\cdot 2^{3h_1+h_2-(3l(h_1)+l(h_2))}
\leq 3\cdot 2^{o(g)-(3l(h_1)+l(h_2)-l(g))} < 2^{o(g)}$ 
($l(g)\leq 2$) so $G$ cannot be optimal. Otherwise, one first remarks that 
$M(H_i)\geq 3$, so $\mu_1(H_i)\leq \frac{1}{2}$ unless 
$h_i=9\cdot 2^{n_i}+u_i$ ($u_i=0,1,2$); but that is impossible for $h_1$ since $l(h_1)=1$. Then one has $|\Aut~G|\leq 3\cdot \frac{1}{2}\cdot 2^{3o(h_1)+o(h_2)}\leq \frac{3}{2}\cdot 2^{o(g)-(3l(h_1)+l(h_2)+1-l(3h_1+h_2+1))}\leq \frac{3}{4}\cdot2^{o(g)}$ (due to (\ref{estimate.l.function})) so $G$ could not have been optimal.
\end{itemize}

If $n_1=1$ and $n_2=2$, then $s_1=2$, $s_2=1$ and $k=2$, then:
\begin{itemize}
 \item If $g=9$ then $2h_1+h_2-3+1+2=9$ so only $h_1=h_2=3$ is possible, i.e. all components are tetrahedra (after stabilization). Then $|\Aut~G|\leq 2\cdot 2^6<|\Aut~C_9|$ so $G$ is not optimal.

 \item If $g=2^m$ or $3\cdot 2^m$ then $2h_1+h_2=g$ and 
$|\Aut~G|\leq 2\cdot 2^{2o(h_1)+o(h_2)}\leq 2^{o(g)-(2l(h_1)+l(h_2)-l(g))}<
2^{o(g)}$ so again $G$ is not optimal.

 \item If $g=3(2^m+2^p)$ with $m\geq p+5$, then the previous estimate
shows that $G$ could only be optimal if $l(h_1)=l(h_2)=1$. Inspecting the
possibilities, one shows that $h_1=3\cdot 2^{m-1}$ and $h_2=3\cdot 2^p$. Now 
the induction shows that $H_2$, of genus $3\cdot 2^p$,
must be a A$_p$ in order to reach equality above (otherwise $\pi(H_2)<1$). 
But then one may compute $|\Aut~G|<2^{o(g)}$ so $G$ is not optimal.

 \item In all other cases, $G$ is made up of a core linked to each of two isomorphic components by an edge. Again pairing the isomorphic components, linking the ends of the edges incident to the core, and reattaching this to one of the original pinching points preserves (or increases) the number of automorphisms, but yields a smaller minimal orbit.
\end{itemize}

We are left then with $n_1=n_2=1$, so $s_1=s_2=k=1$, i.e. two 
(non-isomorphic) components linked by an edge, then 
$|\Aut~G|\leq 2^{o(h_1)+o(h_2)}= 2^{o(g)-(l(h_1)+l(h_2)-l(g))}$.
\begin{itemize}
\item If $g=9$, one can check this configuration cannot be optimal.
\item If $g=2^m$ or $g=3\cdot 2^m$, then $l(h_1)+l(h_2)-l(g))\geq 1$, so
$G$ cannot be optimal.
\item If $g=3(2^m+2^p)$ with $m\geq p+5$, then $2^o(g)$ automorphisms are 
obtained only if
$l(h_1)=l(h_2)=1$. This forces $h_1=3\cdot 2^m$ and $h_2=3\cdot 2^p$ (up 
to swapping $h_1$ and $h_2$). By induction, $G$ may only be $A_m$ and $A_p$
linked by an edge (some checking rules out the case that $p\leq 1$).
\item In all other cases, $l(h_1+h_2)=l(h_1)+l(h_2)$ is forced by the
optimality of $G$. We note that this prohibits $h_1=h_2$; in particular, an optimal graph $G$ cannot be built out of two non-isomorphic components of the same genus, linked by an edge. Thus, if $g=3(2^m+2^p)$, the bound $2^{o(g)}$ may
only be obtained by linking an $A_m$ and an $A_p$ at their roots.
\end{itemize}

{\bf Case 3: $G'$ is disconnected and all its components are isomorphic}. Let $s$ be the number of connected components of $G'$, all of genus $h$, $t$ the incidence degree between two neighbors, and $m$ the number of neighbors of a given component; then $g=sh-s+k+1$ and $2k=smt$; moreover, $s=2$ if and only if $m=1$. We bound, as before (but with the extra possibility of flipping the edges in $O(e)$) $|\Aut~G|\leq 2k\cdot 2^{so(h)}$. The inequality $2k\cdot 2^{so(h)}\leq 2^{o(sh-s+k+1)}$ is equivalent to $2k\leq 2^{o(k+1-s)+sl(h)-l(sh)}=2^{o(k+1-s)+B(s,h)}$.

If $m\geq 3$ (so $s\geq 3$) and $t\geq 2$, then $B(s,h)\geq 1$ by (\ref{long.table}), and $s=\frac{2k}{mt}\leq \frac{k}{3}$; then the inequality is implied by $k\leq 2^{o(\lceil \frac{2k}{3}\rceil+1)}$ which is strict for all $k$ by (\ref{daily.inequalities}).

If $m=2$ and $t\geq 2$, then $k=st\geq 6$, $B(s,h)\geq 1$, 
$s\leq \frac{k}{2}$ and the inequality is implied by 
$k\leq 2^{o(\lceil \frac{k}{2}\rceil +1)}$ which is strict for all $k$ 
except $k=2,4,8$ by (\ref{daily.inequalities}). 
$k=2,4$ are not possible here, while $k=8$ forces $s=4,t=2$. Then $g=4h+5$ 
impossible for $g=9$ or $g=2^n$; in other genera, $B(4,h)\geq 3$ so the 
inequality is easily seen to be strict anyway; no such $G$ may be optimal.

If $t=1$ and $m\geq 4$, then $s\geq 4$ so $B(s,h)\geq 3$ and $s\leq \frac{2k}{m}\leq \frac{k}{2}$ so the inequality is implied by $\frac{k}{4}\leq 2^{o(\lceil \frac{k}{2}\rceil+1)}$ which is strict for all $k$ by (\ref{daily.inequalities}).

If $t=1$ and $m\geq 3$ then $s\geq 4$ so $B(s,h)\geq 3$ and $s\leq \frac{2k}{3}$; we readily get a strict inequality by (\ref{daily.inequalities}).

If $t=1$ and $m=2$ then $G$ is a cycle of components of genus $h$ (each 
component is incident to precisely other two by an edge); $s=k\geq 3$ 
(otherwise $t=2$, in fact), $g=sh+1$. 

\begin{itemize}
\item This is easily seen to be impossible 
for $g\leq 9$ ($h\leq 2$ is forced). 
\item If $g=2^m$, $3\cdot 2^m$, or $g=3(2^m+2^p)$,
$|\Aut~G|\leq 2k\cdot 2^{ko(h)}\leq 2k\cdot 2^{o(g)-(kl(h)+1-l(kh+1))}
\leq 2^{o(g)}\cdot\frac{2k}{2^{B(k,h)}}$. If $k\geq 8$ or $k=6$, the last
fraction is strictly less than one (because at least one of the inequalities
$2^{B(k,h)}\geq 2^{\lceil\frac{k+1}{2}\rceil+1}\geq 2k$ is strict), so $G$ is 
not optimal. For $k=7$, one sees directly that $2^{B(k,h)}>14$ by 
(\ref{long.table}). For $k=4,5$, the inequality $2^{B(k,h)}>2k$ is easily
seen to be true except for $l(h)=1$. However, in that case, the induction
shows that the components have to be graphs of type $A_m$ or $B_p$. In this
case, $|\Aut~G|\leq k\cdot 2^{ko(h)}$, and this is easily seen to be less
than $2^{o(g)}$. For $k=3$ and $l(h)\geq 3$, we have $2^{B(k,h)}>6$. Moreover, $g=3\cdot 2^m$ and $g=3(2^m+2^p)$ cannot be written in the form $3h+1$, so only the case $g=2^m=3h+1$, $k=3$, $l(h)=1$ is left. This is easily seen to be impossible (numerically).

\item In all other genera, we study $2k\leq 2^{B(k,h)}$. $l(h)\geq 3$ would imply strict inequality, so only $l(h)\leq 2$ is possible. $l(h)\geq 2$ and $s\geq 3$ imply by (\ref{main.inequality}) and (\ref{daily.inequalities}) that $G$ is not optimal (split in $s\geq 4$ and $s=3$). So only $l(h)=1$ is left, for which again (\ref{main.inequality}) and (\ref{daily.inequalities}) show that the only hope for optimality is $k=3$. This is only possible for $h=3, g=10$ - a pseudocycle of $K_{3,3}$s with two edges pinched of length three. 
given below.

If the two ends of edges in $O(e)$ incident to a given component are in the 
same orbit, then all ends of edges in $O(e)$ are in the same orbit. These 
ends cannot pinch the same edge of $H$ by (\ref{problems.stabilization}); 
also, (\ref{reduction}) says that we may have at most one ``problem'' when 
trying to stabilize a strictly optimal graph $G$; since we have six incidence 
loci, $H$ is a simple cubic graph. Then we must have 
$M(H)\geq 2$. Then $h\geq 9$ must be one of the special genera of the Main 
Theorem, and for $l(h)=1$ (and $M(H)\geq 2$) we must have 
$\pi(H)\leq \frac{1}{2}2^{o(h)}$.
Then $|\Aut G|\leq 6\cdot \frac{1}{8}2^{3o(h)}=\frac{3}{4}\cdot 2^{o(g)-(3l(h)+1-l(3h+1))}<2^{o(g)}$ so $G$ is not optimal.
If $l(h)=1$ and $h\leq 8$, then $h=3,4,6,8$; we note that for $H$ it must be 
true that fixing one of the two pinched edges is the same as fixing both 
pinched edges; it is easy now to determine that the only genus for which 
there exists a graph $H$ with the maximum number of automorphisms preserving
two edges equal to $\pi(H)=2^{o(h)}$ and $M(H)\geq 2$ is $g=3$, and $H$ must 
be a tetrahedron (pinching this in two opposite edges gives the $K_{3,3}$ with
an edge removed). Thus this last 
situation occurs indeed only for $g=10,h=3$.
\end{itemize}

If $t=1$ and $m=1$ then $G$ is made up of two isomorphic components linked by 
an edge. 
If $g=9$, this is not possible since nine is odd. If $g=8$, $2^{o(g)}$ automorphisms
may only be obtained when the two components of $G'$ 
are tetrahedra each with a pinched edge, i.e. a graph $B_3$ (note that this
graph is not optimal). 
In the other cases, we have
$|\Aut~G|\leq 2\cdot 2^{2o(h)}=2\cdot 
2^{o(g)-(2l(h)-l(g))}\leq 2^{o(g)}$ with equality only if $l(h)=1$.
If 
$g=2^m$ or $3\cdot 2^m$, then the inductive argument yields the unique shape of $G$ as 
$B_m$, respectively $A_m$, proving Theorem B in this case.
If $g=3(2^m+2^p)$, such a graph is not optimal since $l(h)=2$.
\end{proof}

This completes the proof of Theorems A and B.

\end{document}